\documentclass{gtmon_a}
\pdfoutput=1

\usepackage{pinlabel}

\proceedingstitle{Heegaard splittings of 3-manifolds (Haifa 2005)}
\conferencestart{10 July 2005}
\conferenceend{19 July 2005}
\conferencename{Heegaard splittings of 3-manifolds}
\conferencelocation{Haifa}

\editor{Cameron Gordon}
\givenname{Cameron}
\surname{Gordon}

\editor{Yoav Moriah}
\givenname{Yoav}
\surname{Moriah}

\title{On the Heegaard splittings of amalgamated 3--manifolds}
\author{Tao Li}
\givenname{Tao}
\surname{Li}

\address{Department of Mathematics \\
 Boston College \\\newline
 Chestnut Hill, MA 02467\\USA}
\email{taoli@bc.edu}
\urladdr{http://www2.bc.edu/~taoli}

\volumenumber{12}
\issuenumber{}
\publicationyear{2007}
\papernumber{06}
\startpage{157}
\endpage{190}

\doi{}
\MR{}
\Zbl{}

\arxivreference{math.GT/0701395} 


\received{19 November 2006}
\revised{31 August 2007}
\accepted{18 September 2007}
\published{3 December 2007}
\publishedonline{3 December 2007}
\proposed{}
\seconded{}
\corresponding{}
\version{}

\makeop{int}

\makeatletter
\def\cnewtheorem#1[#2]#3{\newtheorem{#1}{#3}[section]
\expandafter\let\csname c@#1\endcsname\c@theorem}

\theoremstyle{plain}
\newtheorem{theorem}{Theorem}[section]
\cnewtheorem{lemma}[theorem]{Lemma}
\cnewtheorem{corollary}[theorem]{Corollary}
\cnewtheorem{proposition}[theorem]{Proposition}
\cnewtheorem{conjecture}[theorem]{Conjecture}
\newtheorem*{theorem1}{\fullref{T02}}
\newtheorem*{theorem*}{Theorem}

\theoremstyle{definition}
\newtheorem{definition}[theorem]{Definition}

\theoremstyle{remark}
\cnewtheorem{remark}[theorem]{Remark}
\cnewtheorem{assumption}[theorem]{Assumption}

\newtheorem*{notation}{Notation}
\makeatother

\numberwithin{figure}{section}

\begin{document}

\begin{htmlabstract}
We give a combinatorial proof of a theorem first proved by Souto which
says the following. Let M<sub>1</sub> and M<sub>2</sub> be simple 3&ndash;manifolds with
connected boundary of genus g>0. If M<sub>1</sub> and M<sub>2</sub> are glued via a
complicated map, then every minimal Heegaard splitting of the
resulting closed 3&ndash;manifold is an amalgamation.  This proof also
provides an algorithm to find a bound on the complexity of the gluing
map.
\end{htmlabstract}

\begin{abstract}
We give a combinatorial proof of a theorem first proved by Souto which
says the following. Let $M_1$ and $M_2$ be simple 3--manifolds with
connected boundary of genus $g>0$. If $M_1$ and $M_2$ are glued via a
complicated map, then every minimal Heegaard splitting of the
resulting closed 3--manifold is an amalgamation.  This proof also
provides an algorithm to find a bound on the complexity of the gluing
map.
\end{abstract}

\begin{asciiabstract}
We give a combinatorial proof of a theorem first proved by Souto which
says the following. Let M_1 and M_2 be simple 3-manifolds with
connected boundary of genus g>0. If M_1 and M_2 are glued via a
complicated map, then every minimal Heegaard splitting of the
resulting closed 3-manifold is an amalgamation.  This proof also
provides an algorithm to find a bound on the complexity of the gluing
map.
\end{asciiabstract}

\maketitle

\section{Introduction}
 
The study of Heegaard splitting has been dramatically changed since
Casson and Gordon introduced the notion of strongly irreducible
Heegaard splitting \cite{CG}. Casson and Gordon proved in \cite{CG}
that if a Heegaard splitting is irreducible but weakly reducible, then
one can perform some compressions on both sides of the Heegaard
surface and obtain an incompressible surface.

Conversely, let $F$ be a connected separating incompressible surface
in a closed 3--manifold $M'$ and $M_1$ and $M_2$ the two manifolds
obtained by cutting open $M'$ along $F$. Then one can construct a
weakly reducible Heegaard splitting by amalgamating two splittings of
$M_1$ and $M_2$ along $F$, see Scharlemann \cite{S2} for more detailed
discussion.
 
In \cite{La} Lackenby showed that if $M_1$ and $M_2$ are simple and
the gluing map is a high power of a pseudo-Ansov homeomorphism of $F$
($F$ is connected), then the minimal genus Heegaard splitting of $M'$
is obtained from splittings of $M_1$ and $M_2$ by amalgamation.  This
implies that the genus of $M'$ is $g(M_1)+g(M_2)-g(F)$.

As pointed out in \cite{La}, it is generally believed that the same is
true if the gluing map is of high distance in the curve complex, see
\fullref{Tmain}.  Note that a high power of a pseudo-Ansov map has
high distance in the curve complex.  Souto \cite{So} proved this first
using the same principles as in \cite{La} by analyzing the geometry
near the incompressible surface.  In this paper, we give a
combinatorial proof of this result and this proof also provides an
algorithm to find the bound on the distance for the gluing map.

\begin{theorem}\label{Tmain}
Let $M_1$ and $M_2$ be orientable simple 3--manifolds with connected
boundary and suppose $\partial M_1\cong\partial M_2\cong F$.  Then
there is a finite set of curves $\mathcal{C}_i\subset\partial M_i$ and
a number $N$ such that, if a homeomorphism $\phi\co \partial
M_1\to\partial M_2$ satisfies
$d_{\mathcal{C}(F)}(\phi(\mathcal{C}_1),\mathcal{C}_2)>N$, where
$d_{\mathcal{C}(F)}$ is the distance in the curve complex
$\mathcal{C}(F)$ of $F$, then
\begin{enumerate}
\item every minimal genus Heegaard splitting of $M_1\cup_\phi M_2$ is
an amalgamation,
\item the Heegaard genus satisfies $g(M_1\cup_\phi
M_2)=g(M_1)+g(M_2)-g(F)$.
\end{enumerate}
Moreover, there is an algorithm to find $\mathcal{C}_i$ and $N$.
\end{theorem}

In this paper, we will study 0--efficient triangulations for
3--manifolds with connected boundary. A 0--efficient triangulation for
a 3--manifold with connected boundary is a triangulation with only one
vertex (on the boundary), the only normal disk is vertex linking, and
there is no normal $S^2$.  By an in-depth analysis of normal annuli in
such triangulations, we prove the following theorem which can be
viewed as a generalization of Hatcher's theorem \cite{Ha} and a
theorem of Jaco and Sedgwick \cite{JS} to manifolds with higher genus
boundary.

\begin{theorem}\label{T02}
Let $M$ be a simple 3--manifold with connected boundary and a
0--efficient triangulation.  Let $S_k$ be the set of normal and almost
normal surfaces satisfying the following two conditions
\begin{enumerate} 
\item the boundary of each surface in $S_k$ consists of essential
curves in $\partial M$
\item the Euler characteristic of each surface in $S_k$ is at least
$-k$.
\end{enumerate} 
Let $C_k$ be the set of boundary curves of surfaces in $S_k$.  Then
$C_k$ has bounded diameter in the curve complex of $\partial M$.
Moreover, there is an algorithm to find the diameter.
\end{theorem}

Since every incompressible and $\partial$--incompressible surface in
$M$ is isotopic to a normal surface in any triangulation, an immediate
corollary of \fullref{T02} is that the set of boundary curves of
essential surfaces with bounded Euler characteristic has bounded
diameter in the curve complex of $\partial M$.

It seems that a version of \fullref{Tmain} is true without the
assumption that $M_i$ is atoroidal or $\partial M_i$ is
incompressible.

\begin{conjecture}\label{Camal}
Let $M_i$ ($i=1,2$) be an irreducible 3--manifold with connected
boundary $\partial M_1\cong\partial M_2\cong F$. Let $\mathcal{D}_i$
be the set of essential curves in $F$ that bound disks in $M_i$.  Then
there is an essential curve $\mathcal{C}_i$ ($i=1,2$) in $\partial
M_i$ such that if the distance between
$\mathcal{D}_2\cup\mathcal{C}_2$ and
$\phi(\mathcal{D}_1\cup\mathcal{C}_1)$ in the curve complex
$\mathcal{C}(F)$ is sufficiently large, then any minimal-genus
Heegaard splitting of $M_1\cup_\phi M_2$ can be constructed from an
amalgamation.
\end{conjecture} 

This conjecture can be viewed as a generalization of
\fullref{Tmain} and a theorem of Scharlemann and Tomova
\cite{ST1}.  Note that in the case that both $M_1$ and $M_2$ are
handlebodies, $\mathcal{C}_1$ and $\mathcal{C}_2$ can be chosen to be
empty and the theorem of Scharlemann and Tomova \cite{ST1} can be
formulated as: if the distance between $\mathcal{D}_2$ and
$\phi(\mathcal{D}_1)$ (ie, the Hempel distance) is large, then the
genus of any other Heegaard splitting must be large unless it is a
stabilized copy of $F$.

The proof in this paper is different from the original proof presented
in the Haifa workshop in 2005, though both proofs use Jaco and
Rubinstein's theory on 0--efficient triangulation \cite{JR}.  This
proof is a byproduct of an effort of finding an algorithmic proof of
the generalized Waldhausen conjecture \cite{L1} and it gives a much
clearer algorithm than the original proof.

I would like to thank Saul Schleimer and Dave Bachman for useful
conversation about their work \cite{BSS} with Eric Sedgwick.  I also
thank the referee for many helpful comments and corrections. The
research was partially supported by an NSF grant.

Throughout this paper, we will denote the interior of $X$ by $\int(X)$,
the closure (under path metric) of $X$ by $\overline{X}$, and the
number of components of $X$ by $|X|$.

\section{Strongly irreducible Heegaard surfaces}\label{Sstr}

Let $M_1$ and $M_2$ be orientable simple 3--manifolds with connected
boundary and suppose $\partial M_1\cong\partial M_2\cong F$.  Let
$\phi\co \partial M_1\to\partial M_2$ be a homeomorphism and $M'=
M_1\cup_\phi M_2$ the closed manifold by gluing $M_1$ and $M_2$ via
$\phi$.  Thus there is an embedded surface $F$ in $M'$ such that
$\overline{M'-F}$ is the disjoint union of $M_1$ and $M_2$.  Since
each $M_i$ is irreducible and $\partial M_i$ is incompressible in
$M_i$, $M'$ is irreducible and $F$ is incompressible in $M'$.  We may
regard $M_i$ as a submanifold of $M'$.

From any Heegaard splittings of $M_1$ and $M_2$, one can naturally
construct a Heegaard splitting of $M'$, called amalgamation.  This
operation was defined by Schultens \cite{Sch}. We give a brief
description below, see \cite{S2,Sch} for details.  Any Heegaard
surface $S_i$ of $M_i$ ($i=1,2$) decomposes $M_i$ into a handlebody
and a compression body.  Each compression body can be obtained by
attaching 1--handles to $F\times I$, a product neighborhood of $F$.
One can extend the 1--handles of the compression body of $M_1$
vertically through the product region $F\times I$ and attach these
extended 1--handles to the handlebody in the splitting of $M_2$.  This
operation produces a handlebody of genus $g(S_1)+g(S_2)-g(F)$.  It is
easy to check that its complement is also a handlebody and we get a
Heegaard splitting of $M'$.  This Heegaard splitting is called an
\emph{amalgamation} of $S_1$ and $S_2$.  Clearly, the resulting
Heegaard splitting from amalgamation is weakly reducible, see
\cite{CG,S2} for definitions and basic properties of weakly reducible
and strongly irreducible Heegaard splittings.

Given a weakly reducible but irreducible Heegaard surface $S$, Casson
and Gordon \cite{CG} showed that one can compress $S$ on both sides
along a maximal collection of disjoint compressing disks and obtain an
incompressible surface.  Scharlemann and Thompson generalized this
construction and gave a construction of \emph{untelescoping} of a
weakly reducible Heegaard splitting, see \cite{S2,ST} for details.
The following lemma follows trivially from the untelescoping
construction.  We say two surfaces intersect nontrivially if they
cannot be made disjoint by an isotopy.

\begin{lemma}\label{Luntel}
Let $S$ be an irreducible Heegaard surface of $M_1\cup_\phi M_2$.
Then either
\begin{enumerate}
\item $S$ is an amalgamation of two splittings of $M_1$ and $M_2$, or
\item there is a submanifold $M_F$ of $M_1\cup_\phi M_2$ ($M_F$ may be
$M_1\cup_\phi M_2$) such that $F\subset M_F$ and $\partial M_F$, if
non-empty, is incompressible, and there is a strongly irreducible
Heegaard surface $S'$ of $M_F$ such that the genus of $S'$ is at most
$g(S)$ and $S'$ nontrivially intersects $F$, or
\item there is an incompressible surface $S'$ with genus less than
$g(S)$ such that $S'$ nontrivially intersects $F$.\qed
\end{enumerate}
\end{lemma}

\section{Intersection with $F$}\label{Sinter}

Suppose $S$ is a minimal genus Heegaard surface of $M'=M_1\cup_\phi
M_2$.  So the genus $g(S)$ is at most $g(M_1)+g(M_2)-g(F)$.  To prove
\fullref{Tmain}, ie, $S$ is an amalgamation, we need to rule out
case (2) and case (3) in \fullref{Luntel}.

We first consider the case (3) in \fullref{Luntel}.  The following
lemma is easy to prove.

\begin{lemma}\label{Lincomp}
Suppose there is an incompressible surface $S'$ that nontrivially
intersects $F$.  Then there is an incompressible and
$\partial$--incompressible surface $S_i$ in $M_i$ such that
$d_{\mathcal{C}(F)}(\phi(\partial S_1),\partial S_2)< -\chi(S')$.
\end{lemma}
\begin{proof}
Since both $S'$ and $F$ are incompressible, we may assume that $S'\cap
F$ consists of essential curves.  Let $S_i'=S'\cap M_i$ ($i=1,2$).
Hence $S_i'$ is incompressible in $M_i$ and
$d_{\mathcal{C}(F)}(\phi(\partial S_1'),\partial S_2')=0$.

If $S_1'$ is $\partial$--compressible, then we perform a
$\partial$--compression on $S_1'$ and get a new incompressible surface
$S_1''$.  Clearly $d_{\mathcal{C}(F)}(\partial S_1',\partial S_1'')\le
1$.  Note that $S_1'$ and $S_1''$ are not $\partial$--parallel in
$M_1$, because otherwise $S'$ can be isotoped to be disjoint from $F$,
contradicting our hypothesis.  Thus, after fewer than $-\chi(S_1')$
$\partial$--compressions, we obtain an incompressible and
$\partial$--incompressible surface $S_1$ in $M_1$. So
$d_{\mathcal{C}(F)}(\partial S_1',\partial S_1)<-\chi(S_1')$.
Similarly, we can find an incompressible and
$\partial$--incompressible surface $S_2$ in $M_2$ with
$d_{\mathcal{C}(F)}(\partial S_2',\partial
S_2)<-\chi(S_2')$. Therefore, $d_{\mathcal{C}(F)}(\phi(\partial
S_1),\partial S_2)< -\chi(S_1')-\chi(S_2')=-\chi(S')$.
\end{proof}

Next we consider the case (2) in \fullref{Luntel}.  Bachman,
Schleimer and Sedgwick \cite{BSS} proved a version of
\fullref{Lincomp} for strongly irreducible Heegaard surfaces, see
\fullref{LBSS} below and \cite[Lemma 3.3]{BSS}.

We first give some definitions using the terminology in \cite{BSS}.

\begin{definition}\label{Dst}
A properly embedded surface is \emph{essential} if it is
incompressible and $\partial$--incompressible.  A properly embedded,
separating surface is {\it strongly irreducible} if there are
compressing disks for it on both sides, and each compressing disk on
one side meets each compressing disk on the other side.  It is {\it
$\partial$--strongly irreducible} if
\begin{enumerate}
	\item every compressing and $\partial$--compressing disk on
	one side meets every compressing and $\partial$--compressing
	disk on the other side, and \item there is at least one
	compressing or $\partial$--compressing disk on each side.
\end{enumerate}
\end{definition}

\begin{lemma}[Bachman--Schleimer--Sedgwick \cite{BSS}]\label{LBSS}
Let $M_F$ be a compact, irreducible, orientable 3--manifold with
$\partial M_F$ incompressible, if non-empty. Suppose $M_F=V \cup _S
W$, where $S$ is a strongly irreducible Heegaard surface. Suppose
further that $M_F$ contains an incompressible, orientable, closed,
non-boundary parallel surface $F$. Then either
\begin{itemize}
	\item $S$ may be isotoped to be transverse to $F$, with every
component of $S - N(F)$ incompressible in the respective submanifold
of $M_F - N(F)$, where $N(F)$ is a small neighborhood of $F$ in $M_F$,
\item $S$ may be isotoped to be transverse to $F$, with every
component of $S - N(F)$ incompressible in the respective submanifold
of $M_F - N(F)$ except for exactly one strongly irreducible component,
or \item $S$ may be isotoped to be almost transverse to $F$ (ie, $S$
is transverse to $F$ except for one saddle point), with every
component of $S - N(F)$ incompressible in the respective submanifold
of $M_F - N(F)$.
\end{itemize}
\end{lemma}

\begin{corollary}\label{CBSS}
Let $M_F$ be a 3--manifold with incompressible boundary and let $F$ be
a separating incompressible and non-boundary parallel surface in
$M_F$.  Let $M_1'$ and $M_2'$ be the 3--manifolds obtained by cutting
$M_F$ open along $F$ and $\phi\co \partial M_1'\to\partial M_2' $ the
gluing map so that $M_F=M_1'\cup_\phi M_2'$.  Suppose $S$ is a
strongly irreducible Heegaard surface of $M_F$. Then there are
surfaces $S_i$ in $M_i'$ such that $d_{\mathcal{C}(F)}(\phi(\partial
S_1),\partial S_2)< -\chi(S)$ and each $S_i$ is either essential or
strongly irreducible and $\partial$--strongly irreducible in $M_i'$.
\end{corollary}
\begin{proof}
Note that if $S\cap M_i'$ consists of $\partial$--parallel surfaces in
$M_i'$ ($i=1,2$), then we can perform an isotopy on $S$ so that $S\cap
F=\emptyset$ after the isotopy, contradicting our hypotheses.  So at
least one component of $S\cap M_i'$ is not $\partial$--parallel.
Moreover, if a component of $S\cap M_i'$ is a disk, since $F$ is
incompressible, the disk must be $\partial$--parallel in $M_i'$ and we
can perform an isotopy on $S$ and remove a trivial-curve component of
$S\cap F$.  Thus, after some isotopies on $S$, we may assume that no
component of $S\cap M_i'$ is a disk. This implies that every curve of
$S\cap F$ is essential in $S$ and hence every component of $S-F$ or
$S\cap M_i'$ is an essential subsurface of $S$.

By \fullref{LBSS}, we can find a component of $S\cap M_i'$, denoted
by $S_i'$ ($i=1,2$), such that (1) each $S_i'$ is an essential
subsurface of $S$ and not $\partial$--parallel in $M_i'$, (2) each
$S_i'$ is either incompressible or strongly irreducible in $M_i'$, and
(3) $d_{\mathcal{C}(F)}(\phi(\partial S_1'),\partial S_2')\le 1$.  If
the third case in \fullref{LBSS} occurs, then both $S_1'$ and $S_2'$
are incompressible. In the other 2 cases, it follows from
\fullref{LBSS} that every curve in $S\cap F$ must be essential in
$F$.  To see this, suppose $\gamma\subset S\cap F$ is an innermost
trivial curve in $F$, then the disk bounded by $\gamma$ in $F$ is a
compressing disk for $S$.  This means that both $S\cap M_1'$ and
$S\cap M_2'$ have compressible components, a contradiction to
\fullref{LBSS}.  Therefore, $\partial S_1'$ and $\partial S_2'$ are
essential in $F$ in any case.

Suppose $S_i'$ is incompressible but $\partial$--compressible in
$M_i'$. As in the proof of \fullref{Lincomp}, after fewer than
$-\chi(S_i')$ $\partial$--compressions, we obtain an essential surface
$S_i$ with $d_{\mathcal{C}(F)}(\partial S_i',\partial
S_i)<-\chi(S_i')$.

Suppose $S_i'$ is strongly irreducible but not $\partial$--strongly
irreducible in $M_i'$.  We say a $\partial$--compressing disk $D$ of
$S_i'$ is \emph{disk-busting} if every compressing disk on the other
side of $S_i'$ intersects $\partial D$.

We first consider the case that $S_i'$ contains a
$\partial$--compressing disk $D$ that is not disk-busting.  So there
is a compressing disk $D'$ on the other side of $S_i'$ with $D\cap
D'=\emptyset$.  Now we perform a $\partial$--compression along $D$ and
get a new surface, which we denote by $S_i''$.  Since $D'\cap
D=\emptyset$, after the isotopy, $D'$ remains a compressing disk of
$S_i''$.  Note that since $S_i'$ is strongly irreducible, by
definition, there is a compressing disk of $S_i'$ on the same side as
$D$, in fact, a simple cutting-and-pasting argument can show that
there is a compressing disk on the same side as $D$ and disjoint from
$D$.  This means that $S_i''$ is still strongly irreducible.

After a finite number of such $\partial$--compressions, we may assume
every $\partial$--compressing disk of $S_i''$ is disk-busting.  If
$S_i''$ is not $\partial$--strongly irreducible, then there must be a
pair of disjoint $\partial$--compressing disks $D$ and $D'$ on
different sides of $S_i''$.  Since $D\cap D'=\emptyset$, we can
perform $\partial$--compressions along $D$ and $D'$ simultaneously.
Since both $D$ and $D'$ are disk-busting, the resulting surface after
$\partial$--compressions along $D$ and $D'$ is incompressible.

Therefore, after fewer than $-\chi(S_i')$ $\partial$--compressions, we
obtain a surface $S_i$ in $M_i'$ ($i=1,2$) such that each $S_i$ is
either essential or strongly irreducible and $\partial$--strongly
irreducible in $M_i'$ and $d_{\mathcal{C}(F)}(\partial S_i',\partial
S_i)<-\chi(S_i')$.  Similar to the proof of \fullref{Lincomp}, we
have $d_{\mathcal{C}(F)}(\phi(\partial S_1),\partial S_2)< -\chi(S)$.
\end{proof}

The next corollary follows trivially from \fullref{Luntel},
\fullref{Lincomp} and \fullref{CBSS}.

\begin{corollary}\label{CStr}
Let $S$ be an irreducible Heegaard surface of $M_1\cup_\phi
M_2$. Suppose $S$ is not a amalgamation of two splittings of $M_1$ and
$M_2$.  Then there is a properly embedded surface with boundary $S_i$
in $M_i$ such that $d_{\mathcal{C}(F)}(\phi(\partial S_1),\partial
S_2)< -\chi(S)$ and each $S_i$ is either essential or strongly
irreducible and $\partial$--strongly irreducible in $M_i$.\qed
\end{corollary}

\begin{remark}\label{R1}
Let $S_1$ and $S_2$ be components of $S\cap M_i$ as in
\fullref{CStr}.  It follows from the construction above and
\fullref{Dst} that the boundary of $S_1$ and $S_2$ consists of
essential curves.  We fix a 0--efficient triangulation (described
below) for each $M_i$. If $S_i$ is essential, then $S_i$ is isotopic
to a normal surface.  If $S_i$ is $\partial$--strongly irreducible,
then by a theorem of Bachman \cite{B}, $S_i$ is isotopic to a normal
or an almost normal surface with boundary.  The referee pointed out a
controversy in a theorem in \cite{B}.  In our proof, we will use
Bachman's theorem, but give a workaround in the appendix avoiding the
controversial part of Bachman's argument.  If $S_i$ is
$\partial$--strongly irreducible, the general case follows from the
appendix is that, after isotopy or possible $\partial$--compressions,
(1) $S_i$ is normal or almost normal and $\partial S_i$ consists of
normal curves in $\partial M_i$, and (2) at most one component of
$\partial S_i$ is a trivial curve and at least one component of
$\partial S_i$ is an essential curve.  Note that a trivial normal
curve in a one-vertex triangulation of $\partial M_i$ is
vertex-linking, see \fullref{P1} below.  For simplicity, we
will assume that $S_i$ is almost normal and $\partial S_i$ consists of
essential normal curves and the proof for the general case is
basically the same.
\end{remark}

\section{The 0--efficient triangulation}\label{S0-eff}

Let $S$ be a minimal genus Heegaard surface.  By \fullref{CStr},
if a Heegaard surface $S$ is not obtained from amalgamation, then
there is a surface $S_i$ properly embedded in $M_i$ such that $S_i$ is
either essential or $\partial$--strongly irreducible in $M_i$ and
$d_{\mathcal{C}(F)}(\phi(\partial S_1),\partial S_2)< -\chi(S)\le
2(g(M_1)+g(M_2)-g(F))-2$.  Fixing a 0--efficient triangulation
(described below) of $M_i$, as in \fullref{R1}, we may assume $S_i$
is a normal or an almost normal surface with respect to the
triangulation and $\partial S_i$ consists of essential normal curves
in $\partial M_i$.  Our goal is to prove that the boundary curves of
such (almost) normal surfaces have bounded diameter in the curve
complex of $F=\partial M_i$, see \fullref{T02}.

The 0--efficient triangulation, introduced by Jaco and Rubinstein
\cite{JR}, is a very convienient tool, see for example \cite{L1}.  In
this paper we are mainly interested in 0--efficient triangulation for
manifolds with connected and incompressible boundary.  We first give
an overview of the definition and special properties of such a
triangulation.

Since $\partial M_i$ is connected and incompressible in $M_i$, by
\cite{JR}, $M_i$ admits a special triangulation with the following
properties:
\begin{enumerate}
\item the triangulation has only one vertex which lies in $\partial
M_i$
\item the only normal disk is the vertex-linking one,
\item there is no normal $S^2$ in $M_i$
\end{enumerate}
We call such a triangulation a \emph{0--efficient triangulation} for
$M_i$.  It is also shown in \cite{JR} that there is an algorithm to
find such a triangulation.

Similar to 0--efficient triangulations for closed 3--manifolds, such
triangulations have some remarkable properties.  The following lemma
was proved by Jaco and Rubinstein and the proof is basically the same
as the closed case, also see \cite[Lemma 5.1]{L1}.  The proof of the
\fullref{L01} uses a technique in \cite{JR} called \emph{barrier}.
A barrier is basically a 2--complex barrier for the normalization
operations.  We refer the reader to \cite[section 5]{L1} for a brief
explanation and \cite[section 3.1]{JR} for more details.  The proof of
\fullref{L01} is similar in spirit to that of \cite[Lemma 5.1]{L1}.

\begin{lemma}\label{L01}
Let $M_i$ be a simple 3--manifold with connected boundary and
$\mathcal{T}$ a 0--efficient triangulation.  Then every properly
embedded normal annulus with respect to $\mathcal{T}$ is
$\partial$--parallel and incompressible.
\end{lemma}
\begin{remark}\label{RBSS}
$M_i$ does not contain any normal M\"{o}bius band, since the boundary
of a neighborhood of a normal M\"{o}bius band is a normal annulus,
which contradicts the above lemma and the assumption that $M_i$ is
simple.
\end{remark}
\begin{proof}
Let $A$ be a properly embedded normal annulus in $M_i$.  Since $M_i$
is simple, every incompressible annulus is $\partial$--parallel. So it
suffices to prove that $A$ is incompressible.

Suppose $A$ is compressible, then $\partial A$ must be trivial curves
in $\partial M_i$ since $\partial M_i$ is incompressible in $M_i$.
Note that the induced triangulation of $\partial M_i$ has only one
vertex. The only trivial normal curve in a one-vertex triangulation of
$\partial M_i$ is vertex-linking (see part (a) of
\fullref{P1}). Hence $\partial A$ is a pair of parallel
vertex-linking curves.  Let $\gamma_1$ and $\gamma_2$ be the two
components of $\partial A$ and $D_j$ ($j=1,2$) the disk bounded by
$\gamma_j$ in $\partial M_i$.  As $\gamma_1$ and $\gamma_2$ are
parallel, we may suppose $D_1\subset D_2$.

Note that the disk $A\cup D_1$ may not be normal, but $A$ is a barrier
for the normalization operations that make $A\cup D_1$ normal.  So we
can normalize $A\cup D_1$ to a normal disk $\Delta$.  Since the
triangulation is 0--efficient, $\Delta$ is a vertex-linking disk.
Since $A$ is a barrier for the normalizing operations, $A$ must lie in
the 3--ball bounded by $\Delta$ and a disk of $\partial M_i$.
However, there is no normal annulus in a small neighborhood of the
vertex, a contradiction.
\end{proof}

\begin{notation}\label{notation}
To simplify notation, in the remaining of the paper, we use $M$ to
denote either $M_1$ or $M_2$ and $F=\partial M_i$.  Unless specified,
we use $S$ to denote the surface $S_i$ in \fullref{CStr}.  We
fix a 0--efficient triangulation of $M$ and assume $S$ is a normal or
an almost normal surface in $M$ with respect to the 0--efficient
triangulation and $\partial S$ consists of essential normal curves in
$\partial M$.
\end{notation}

Now we consider all the properly embedded normal and almost normal
surfaces in $M$ whose boundary consists of essential curves in
$\partial M$.  Similar to \cite{FO,L1}, there is a finite collection
of branched surfaces each obtained by gluing normal disks and at most
one almost normal piece, such that each of these normal or almost
normal surfaces is fully carried by a branched surface in this
collection.  Moreover, similar to \cite{AL,L1}, since there is no
normal $S^2$ and the only normal disk in this triangulation is
vertex-linking, after taking sub-branched surfaces if necessary, we
may assume no branched surface in this collection carries a normal
disk or normal $S^2$.

Let $B$ be a branched surface in this collection that fully carries
$S$.  So $\partial B$ is a train track in $\partial M$.  We call a
train track a \emph{normal train track} if every curve carried by the
train track is normal with respect to the induced triangulation of
$\partial M$. By the construction, $\partial B$ is a normal train
track.

Note that in the general case from the appendix, $\partial S$ may
contain a single trivial curve, though at least one component of
$\partial S$ is essential.  In this case, we may split $B$ so that a
component of $\partial B$ is an isolated trivial circle and each other
component of $\partial B$ fully carries an essential curve (as
required by part (c) of \fullref{P1}).  For simplicity, as
mentioned in \fullref{R1}, we assume $\partial S$ is essential and
$\partial B$ fully carries $\partial S$.

\eject 
\begin{proposition}\label{P1}\ 
\begin{enumerate}
\item[\rm(a)] A normal simple closed curve in $\partial M$ is trivial if
and only if it is vertex-linking.
\item[\rm(b)] At most one component of $\partial M-\partial B$ is a
monogon.
\item[\rm(c)] The train track $\partial B$ does not carry any trivial
curve.
\end{enumerate}
\end{proposition}
\begin{proof}
Part (a) follows from the fact that the induced triangulation of
$\partial M$ has only one vertex.  To see this, let $\gamma$ be a
normal trivial curve and $D$ the disk bounded by $\gamma$ in $\partial
M$.  Let $e$ be any edge (or 1--simplex) in the induced triangulation
of $\partial M$.  If a component $\alpha$ of $e\cap D$ is an arc in
$\int(e)$, then $\alpha$ is properly embedded in $D$ and cuts $D$ into
two subdisks $D_1$ and $D_2$.  As there is only one vertex, at least
one subdisk, say $D_1$, does not contain the vertex.  Hence the
intersection of $D_1$ and the 1--skeleton of the triangulation
consists of arcs properly embedded in $D_1$.  These arcs cut $D_1$
into subdisks and an outermost subdisk is a bigon with one edge in
$\partial D$ and the other edge in the 1--skeleton.  This means that
$\gamma=\partial D$ is not a normal curve, a contradiction.
Therefore, every component of $e\cap D$ is an arc with one endpoint
the vertex of the triangulation and the other endpoint in $\partial
D$.  This implies that $\gamma=\partial D$ is vertex-linking.

The proof of part (b) is similar.  Since every curve carried by
$\partial B$ is a normal curve, the argument above implies that each
monogon component of $\partial M-\partial B$ must contain the vertex
of the triangulation.  Part (b) follows from that assumption that
there is only one vertex in the triangulation.

Part (c) follows from the assumption that $B$ fully carries $S$ and
$\partial S$ consists of essential curves.  Let $N(\partial B)$ be a
fibered neighborhood of the train track $\partial B$ in $\partial
M$. We may assume $\partial S$ lies in $N(\partial B)$ and is
transverse to the interval fibers of $N(\partial B)$.  Since $\partial
B$ fully carries $\partial S$, after some isotopy and taking multiple
copies of $\partial S$ if necessary, we may assume that the horizontal
boundary of $N(\partial B)$ lies in $\partial S$.  Since each
component of $\partial S$ is essential, this means that no horizontal
boundary component of $N(\partial B)$ is a trivial circle.  In other
words, no component of $\overline{\partial M-\partial B}$ (or
$\partial M-\int(N(\partial B)$) is a disk with smooth boundary.  If
$\partial B$ carries a trivial circle $\gamma$, then a trivial index
argument implies that the disk bounded by $\gamma$ contains either a
disk component of $\overline{\partial M-\partial B}$ with smooth
boundary or at least two monogons.  The first case is impossible by
the argument above and the second case is ruled out by part (b).  So
$\partial B$ does not carry any trivial curve.
\end{proof}

Each surface carried by $B$ is corresponding to a nonnegative integer
solution to the system of branch equations, see \cite{AL,FO,L1} for
more detailed discussion.  To simplify notation, we will not
distinguish between a surface carried by $B$ and its corresponding
integer solution to the system of branch equations.

By the normal surface theory, there is a finite set of fundamental
solutions of the system of branch equations such that any surface
carried by $B$ is a linear combination of the fundamental solutions
with nonnegative integer coefficients.  We denote the fundamental
solutions by $F_1,\dots,F_s, C_1,\dots, C_t, A_1,\dots A_n$, where
each $A_j$ is a normal annulus carried by $B$, each $C_j$ is a closed
surface carried by $B$ and the $F_j$'s are the other fundamental
solutions.  So the surface $S$ can be written as $S=\sum s_jF_j+\sum
t_jC_j +\sum n_jA_j$ where each $s_j$, $t_j$ or $n_j$ is a nonnegative
integer.

\begin{proposition}\label{Pbounded}
$\sum s_j\le 2-\chi(S)$.
\end{proposition}
\begin{proof}
Since $S$ is a normal or an almost normal surface, we may assume that
at most one fundamental solution contains an almost normal piece and
its coefficient in the linear combination above is either 0 or 1.

Note that $M$ does not contain any normal projective plane, since the
boundary of a twisted $I$--bundle over a normal $P^2$ is a normal
$S^2$ and $M$ does not contain any normal $S^2$. Moreover $B$ does not
carry any normal disk by our assumption.  These imply that $B$ does
not carry any normal surface with positive Euler characteristic.

We first consider the case that $S$ is a normal surface.  First, we
have $\chi(S)=\sum s_j\chi(F_j)$ $+\sum t_j\chi(C_j) +\sum n_j\chi(A_j)$.
Since $S$ is normal, each fundamental solution with positive
coefficient in the linear combination above is a normal surface.
Since $B$ does not carry a normal surface with positive Euler
characteristic, we have
$$\chi(S)=\sum s_j\chi(F_j)+\sum t_j\chi(C_j)\le \sum s_j\chi(F_j)\le
-\sum s_j.$$
So in the case that $S$ is a normal surface, we have $\sum s_j\le
-\chi(S)$.  If $S$ is almost normal, we may suppose some $C_k$ (or
$F_k$) is almost normal and the coefficient of $C_k$ (or $F_k$) is 1.
Note that since $\chi(C_k)$ (or $\chi(F_k)$) is at most 2, we have
$\chi(S)\le 2-\sum s_j$ and $\sum s_j\le 2-\chi(S)$.
\end{proof}

Since there are only finitely many such branched surfaces $B$, to
prove \fullref{T02}, it suffices to show that the set of boundary
curves of surfaces carried by $B$ with bounded Euler characteristic
has bounded diameter in the curve complex of $F$.  Since each $C_j$ is
a closed surface, $\partial S=\sum s_j\partial F_j+\sum n_j\partial
A_j$.  As $\sum s_j$ is bounded by \fullref{Pbounded}, there
are only finitely many possibilities for curves $\sum s_j\partial
F_j$.  Thus the key part of the proof is to study normal annuli
carried by $B$.

\section{Normal annuli}\label{Sannuli}

We use the same notation.  Let $B$ be a branched surface in $M$ that
fully carries $S$ as above and $A_1,\dots,A_n$ the fundamental
solutions that correspond to normal annuli carried by $B$.  Since $B$
does not carry any normal surface of positive Euler characteristic,
each component of the normal sum $\sum n_iA_i$ must have Euler
characteristic 0 and hence is either a normal torus or a normal
annulus carried by $B$.  Note that there is no normal Klein bottle in
the 0--efficient triangulation, see Lemma 5.1 and Corollary 5.2 in
\cite{L1}.

Let $N(B)$ be a fibered neighborhood of $B$ and $\pi\co  N(B)\to B$ the
map collapsing each $I$--fiber to a point, see \cite{FO,L1} for more
details. We may view $A_1,\dots A_n$ as embedded annuli in $N(B)$.
Then $\pi(\sum n_iA_i)$ is a sub-branched surface of $B$ fully
carrying $\sum n_iA_i$.  Since each $A_j$ is $\partial$--parallel,
there is an annulus $\Gamma_j\subset\partial M$ such that
$\partial\Gamma_j=\partial A_j$ and $A_j$ is isotopic to $\Gamma_j$
relative to $\partial A_j$.  Throughout this paper, we will use $T_j$
to denote the solid torus bounded by $A_j\cup\Gamma_j$.

Next we study the intersection of two normal annuli carried by $B$.
Let $A_1$ and $A_2$ be two annuli carried by $B$ and suppose $A_1\cap
A_2\ne\emptyset$.  As above, let $\Gamma_1$ and $\Gamma_2$ be the
annuli in $\partial M$ bounded by $\partial A_1$ and $\partial A_2$
respectively.

If $A_1\cap A_2$ contains a closed curve $\gamma$, then since every
normal annulus is incompressible by \fullref{L01}, $\gamma$ is
either trivial in both $A_1$ and $A_2$ or essential in both $A_1$ and
$A_2$.  Let $\Gamma$ be the union of closed curves in $A_1\cap A_2$
that are trivial in both $A_1$ and $A_2$.  Let $P_i$ be the component
of $A_i-\Gamma$ that contains $\partial A_i$ ($i=1,2$).  Clearly
$P_1\cap P_2=(A_1\cap A_2)-\Gamma$.  Now we perform standard cutting
and pasting along $\Gamma$ and denote by $A_i'$ the resulting
component that contains $P_i$ ($i=1,2$).  If $A_1'=A_2'$, then
$\chi(A_1')<0$, which means that the cutting and pasting above also
produces an embedded normal surface with positive Euler
characteristic, a contradiction to the assumptions of the branched
surface $B$.  Thus $A_1'\ne A_2'$, each $A_i'$ is an embedded normal
annulus, and $A_1'\cap A_2'$ does not contain any trivial closed
curves.  Note that the cutting and pasting above may produce a normal
torus.  Therefore, after some cutting and pasting above, we may assume
the intersection of two normal annuli does not contain trivial curves.

\begin{figure}[ht!]
\begin{center}
\labellist\small \pinlabel $\alpha$ [l] at 190 539 \pinlabel $\alpha$
[l] at 469 542 \hair 7pt\pinlabel {type $I$} [t] at 197 410 \pinlabel
{type $II$} [t] at 478 410 \pinlabel (a) [t] <0pt,-15pt> at 197 410
\pinlabel (b) [t] <0pt,-15pt> at 478 410 \endlabellist
\includegraphics[width=4in]{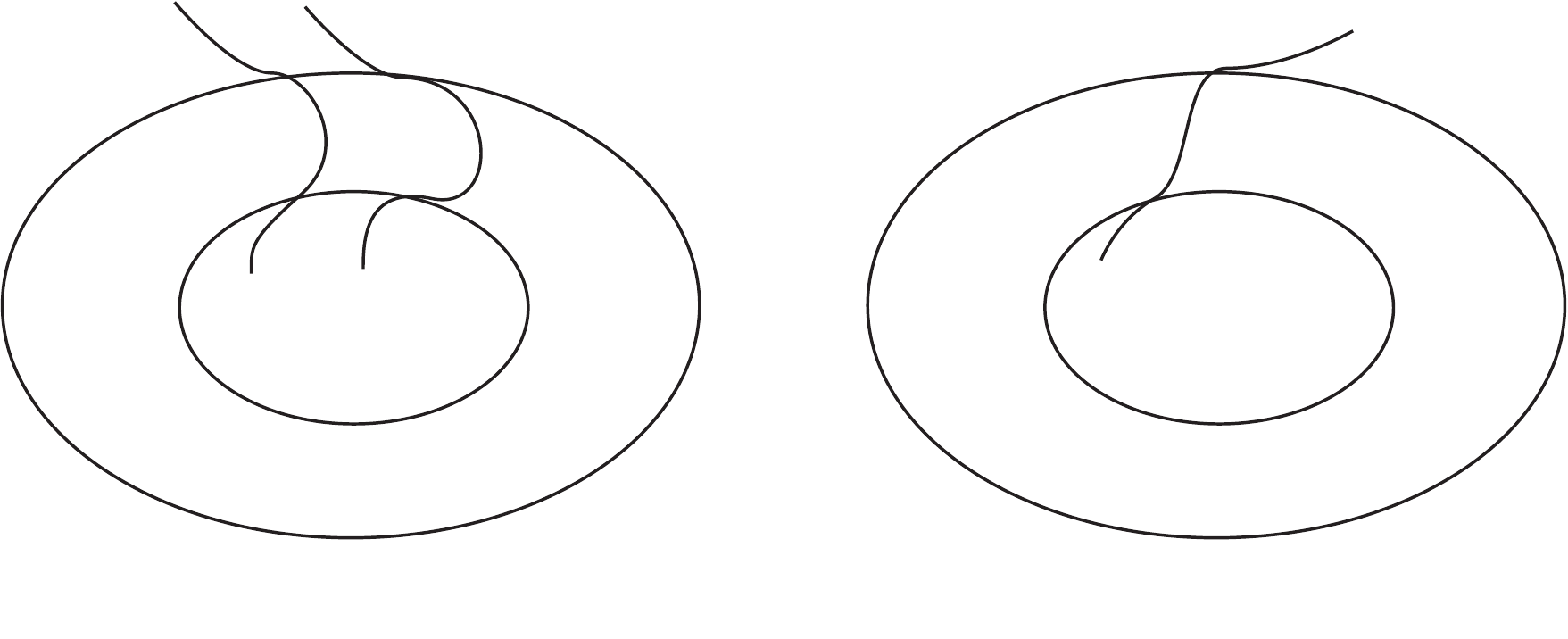}
\caption{}\label{Ftype}
\end{center}
\end{figure}

\begin{definition}\label{Dtype}
Suppose $\partial A_1\cap\partial A_2\ne\emptyset$. This means that
$\Gamma_1\cap\partial A_2\ne\emptyset$.  We consider an arc $\alpha$
of $\Gamma_1\cap\partial A_2$ with endpoints in different components
of $\partial\Gamma_1$.  Since $\partial A_1$ and $\partial A_2$ are
carried by $B$, $\partial A_1\cup\partial A_2$ naturally forms a train
track.  We say $\alpha$ is of type $I$ in $\Gamma_1$ if
$\partial\Gamma_1\cup\alpha$ form a train track of a Reeb annulus, as
shown in \fullref{Ftype}(a).  Otherwise, the train track $\partial
A_1\cup\alpha$ is as shown in \fullref{Ftype}(b) and we say
$\alpha$ is of type $II$ in $\Gamma_1$.  We say $A_1$ is of type $I$
relative to $A_2$ if a component of $\Gamma_1\cap\partial A_2$ is of
type $I$, otherwise we say $A_1$ is of type $II$ relative to $A_2$.
\end{definition}

 Note that if there are two type $I$ arcs of $\Gamma_1\cap \partial
 A_2$ with opposite switching directions along $\partial\Gamma_1$,
 then the train track $\pi(\partial A_1\cup\partial A_2)$ carries a
 trivial circle.  By part (c) of \fullref{P1}, $\partial B$
 does not carry any trivial circle.  So all the type $I$ arcs of
 $\Gamma_1\cap \partial A_2$ must have coherent switching directions
 as shown in \fullref{Ftype}(a), ie, the train track formed by
 $\partial A_1$ and these type $I$ arcs carries a Reeb lamination of
 an annulus.

\begin{proposition}\label{PReeb}
Let $\Gamma$ be an annulus in $\partial M$ and suppose
$\partial\Gamma$ consists of normal curves with respect to a
one-vertex triangulation of $\partial M$. Let $\alpha$ be a properly
embedded essential arc in $\Gamma$. Suppose $\partial\Gamma\cup\alpha$
forms a Reeb train track as shown in \fullref{Ftype}(a) which
carries a normal Reeb lamination.  Then $\Gamma$ contains the vertex
of the triangulation.
\end{proposition}
\begin{proof}
We may deform $\partial\Gamma\cup\alpha$ into a train track $\tau$.
Note that our hypothesis says that the Reeb lamination carried by
$\tau$ is normal with respect to the one-vertex triangulation of
$\partial M$.

Suppose that $\Gamma$ does not contain the vertex.  Let $e$ be an edge
intersecting $\Gamma$.  Let $\beta$ be a component of $e\cap\Gamma$.
Since $\Gamma$ does not contain the vertex, there are only two
possibilities: (1) $\partial\beta$ lies in the same circle of
$\partial\Gamma$ and (2) the endpoints of $\beta$ lie in different
circles of $\partial\Gamma$.  If $\partial\beta$ lies in the same
component of $\partial\Gamma$, then similar to the proof of
\fullref{P1}, this component of $\partial\Gamma$ must have
trivial intersection with the edge $e$ and hence cannot be a normal
curve.  Similarly, if the endpoints of $\beta$ lie in different
circles of $\partial\Gamma$, then every non-compact leaf of the Reeb
lamination has trivial intersection with the edge $e$ and hence the
Reeb lamination carried by $\tau$ cannot be normal.  Thus $\Gamma$
must contain the vertex.
\end{proof}

An isotopy is called a \emph{normal isotopy} if it is invariant on
each simplex of the triangulation.  Next we will perform some normal
isotopies on $\partial A_i$.  If $\gamma$ is normally isotopic to
$\partial A_i$, then $A_i$ is normally isotopic to a normal annulus
$A_i'$ with $\partial A_i'=\gamma$.  Moreover, for any surface $X$
carried by $B$, we may assume $A_i+X$ is normally isotopic to
$A_i'+X$.  Next we will perform some normal isotopies and these normal
isotopies do not change the surface under normal sum.

\begin{definition}
Let $X$ be a point of $\partial A_1\cap\partial A_2$.  A small
neighborhood of $X$ is cut into 4 corners by $\partial A_1\cup\partial
A_2$.  A corner is called a cusp if it becomes a cusp after deforming
$\partial A_1\cup\partial A_2$ into a train track.  We call a disk $D$
in $\partial M$ a \emph{bigon} if (1) $\partial D$ consists of two
arcs, one from $\partial A_1$ and the other from $\partial A_2$, and
(2) the two corners of $D$ at $\partial A_1\cap\partial A_2$ are both
cusps.  We say $D$ is an innermost bigon if $\int(D)\cap (\partial
A_1\cup\partial A_2)=\emptyset$.  A bigon is said to be \emph{trivial}
if it does not contain the vertex of the triangulation.
\end{definition}

\vspace{5pt}

\noindent\textbf{Eliminate a trivial bigon}\qua Let $D$ be an innermost
trivial bigon and $b_1\subset\partial A_1$ and $b_2\subset\partial
A_2$ the two edges of $\partial D=b_1\cup b_2$.  Since $D$ does not
contain the vertex and both $\partial A_1$ and $\partial A_2$ are
normal curves, the intersection of $D$ and the 1--skeleton of the
triangulation consists of arcs with one endpoint in $b_1$ and the
other endpoint in $b_2$. This means that $b_1$ and $b_2$ are normally
isotopic.  Hence we can perform a normal isotopy on $A_i$ near
$\partial A_i$, changing $\partial A_1$ to $(\partial A_1-b_1)\cup
b_2$ and $\partial A_2$ to $(\partial A_2-b_2)\cup b_1$.  After the
normal isotopy and a small perturbation, $\partial A_1\cap\partial
A_2$ has fewer intersection points.  We may successively eliminate all
the trivial bigons using such normal isotopies.

For a given finite set of annuli carried by $B$, after some normal
isotopies as above, we may assume that for any pair $A_i$ and $A_j$,
$\partial A_i\cup\partial A_j$ does not form any trivial bigon.

\begin{definition}
Let $\alpha$ be an arc component of $A_1\cap A_2$ that is trivial (ie
$\partial$--parallel) in both $A_1$ and $A_2$.  Then $\alpha$ together
with a subarc $\beta_i$ of $\partial A_i$ ($i=1,2$) bounds a subdisk
$D_i$ of $A_i$.  If $D_1\cap D_2=\alpha$ then $\beta_1\cup\beta_2$
bounds a disk $\Delta$ in $\partial M$ and $D_1\cup D_2\cup\Delta$ is
a 2--sphere bounding a 3--ball.  We call such a 3--ball a
\textit{football region}.  Note that since the endpoints of $\beta_i$
are also the endpoints of $\alpha$ and since $A_1$ and $A_2$ are
carried by the same branched surface $B$, after deforming
$\beta_1\cup\beta_2$ into train track, $\beta_1\cup\beta_2$ cannot
form a monogon.  Since the train track $\partial B$ does not carry any
trivial circle, $\Delta$ must be a bigon.  Moreover, since we have
assumed that there is no trivial bigon, the bigon $\Delta$ must
contain the vertex of the triangulation.  A football region is said to
be \emph{innermost} if it does not contains any other football region.
A football region bounded by $D_1\cup D_2\cup\Delta$ said to be
\emph{trivial} if $D_1\cap A_2=D_2\cap A_1=\alpha$. Clearly a trivial
football region must be innermost.
\end{definition}

\vspace{5pt}

\noindent\textbf{Eliminate a trivial football region}\qua Suppose the
football region bounded by $D_1\cup D_2\cup\Delta$ as above is
trivial.  Let $\alpha=\partial D_1\cap\partial D_2$. Since $D_1\cap
A_2=D_2\cap A_1=\alpha$, we can perform a canonical cutting and
pasting along $\alpha$ and obtain annuli $(A_1-D_1)\cup D_2$ and
$(A_2-D_2)\cup D_1$. Clearly $(A_1-D_1)\cup D_2$ and $(A_2-D_2)\cup
D_1$ are embedded annuli carried by $B$ and are isotopic to $A_1$ and
$A_2$ respectively.  After a slight perturbation, the resulting annuli
have fewer intersection curves.  Thus, after a finite number of such
operations, we may assume there is no trivial football region.

\begin{definition}
We say $A_1\cup A_2$ is bigon-efficient if $A_1\cap A_2$ contains no
trivial closed curve, $\partial A_1\cup\partial A_2$ does not form any
trivial bigon in $\partial M$, and $A_1\cup A_2$ does not form any
trivial football region.
\end{definition}

As above, we can perform some canonical cutting and pasting along
$A_1\cap A_2$ and get a pair of new annuli $A_1'$ and $A_2'$ such that
$A_1'\cup A_2'$ is bigon-efficient.  By our construction, $A_1'$ and
$A_2'$ are also carried by $B$ and $A_1+A_2=A_1'+A_2'$.

Next we will assume that $A_1\cup A_2$ is bigon-efficient and consider
the intersection pattern of $A_1\cap A_2$.

\begin{lemma}\label{Louter}
Let $\beta_0$ be an arc in $A_1\cap A_2$ and suppose $\beta_0$ is
$\partial$--parallel in $A_1$.  Let $\Delta_0$ be the subdisk of $A_1$
bounded by $\beta_0$ and a subarc of $\partial A_1$. Let
$\beta_1,\dots,\beta_k$ be the components of $\int(\Delta_0)\cap A_2$.
Suppose each $\beta_i$ ($i\ge 1$) is outermost in $A_1$.  Then at
least one $\beta_i$ ($i\ge 1$) is $\partial$--parallel in $A_2$.
\end{lemma}
\begin{proof}
Suppose each $\beta_i$ ($i\ge 1$) is an essential arc in $A_2$.  Let
$\delta_i$ be the subdisk of $\Delta_0$ bounded by $\beta_i$ ($i\ge
1$) and a subarc of $\partial A_1$.  Since $\beta_i$ ($i\ge 1$) is
essential in $A_2$ and outermost in $A_1$, each $\delta_i$ is a
$\partial$--compressing disk for $A_2$.  This implies that
$\partial\delta_i\cap\partial M$ is a type $I$ arc in $\Gamma_2$.  By
\fullref{PReeb}, $\Gamma_2$ contains the vertex of the
triangulation.

Since $A_2$ is $\partial$--parallel in $M$, $A_2\cup\Gamma_2$ bounds a
solid torus $T_2$.  Let $M'$ be the closure of $M-T_2$.  So $M'\cong
M$ and we may view $A_2$ as an annulus in $\partial M'$.

We use $D$ to denote the closure of $\Delta_0-\cup_{i=1}^k\delta_i$.
Thus we may view $D$ as a disk properly embedded in $M'$. Since
$\partial M'$ is incompressible in $M'$, $\partial D$ bounds a disk
$D'$ in $\partial M'$.  We view $A_2$ as a subannulus of $\partial
M'$. So $D'\cap A_2\ne\emptyset$.

Note that $\partial A_2$ cuts $D'$ into disks and at least two such
disks are outermost in $D'$ (an outermost disk is a disk whose
boundary consists of a subarc of $\partial D'$ and a subarc of
$\partial A_2$).  Let $\Delta$ be such an outermost disk.  If
$\Delta\subset A_2\subset\partial M'$, then since each $\beta_i$
($i\ge 1$) is essential in $A_2$, $\beta_0$ must be an arc in
$\partial\Delta$.  Since there are at least two outermost disks, we
may choose $\Delta$ to be outside $A_2$.  In other words,
$\Delta\subset\partial M'-\int(A_2)=\partial M-\int(\Gamma_2)$.  Since
$\Gamma_2$ contains the vertex of the triangulation, this means that
$\Delta$ does not contain the vertex.  If we deform $\partial
A_1\cup\partial A_2$ into a train track, then $\Delta$ becomes either
a bigon or a monogon or a smooth disk.  As in the proof of
\fullref{P1}, a monogon or a smooth disk must contain the
vertex.  Since $\Delta$ does not contain the vertex, $\Delta$ must be
a trivial bigon, which contradicts our assumption that $A_1\cup A_2$
is bigon-efficient.
\end{proof}

\begin{lemma}\label{LFB}
Let $A_1$ and $A_2$ be as above and suppose $ A_1\cup A_2$ is
bigon-efficient.  Then $A_1$ and $A_2$ do not form any football
region.
\end{lemma}
\begin{proof}
Suppose there is a football region $X$ bounded by $D_1\cup
D_2\cup\Delta$, where $D_i\subset A_i$ is a disk bounded by a
component $\alpha$ of $A_1\cap A_2$ and a subarc of $\partial A_i$ and
$\Delta\subset\partial M$.  We use $\beta_i$ ($\beta_i\subset\partial
A_i$) to denote $\partial D_i-\int(\alpha)$ ($i=1,2$).  Note that
$\Delta$ must contain the vertex of the triangulation, because
otherwise $\Delta$ is a trivial bigon contradicting that $A_1\cup A_2$
is bigon-efficient. Without loss of generality, we may assume $X$ does
not contain any other football region.

If $D_1\cap A_2=D_2\cap A_1=\alpha$, then the 3--ball bounded by
$D_1\cup D_2\cup\Delta$ is a trivial football region, contradicting
the assumption that $A_1\cup A_2$ is bigon-efficient. So we may assume
$\int(D_1)\cap A_2\ne\emptyset$.

Since $\int(D_1)\cap A_2\ne\emptyset$, we can always find a component
$\beta_0$ of $D_1\cap A_2$ such that $\beta_0$ is not outermost in
$A_1$ but every component of $\int(D_1)\cap A_2$ inside the disk
bounded by $\beta_0$ and a subarc of $\partial A_1$ is outermost in
$A_1$.  By \fullref{Louter}, there is at least one arc
$\alpha'\subset \int(D_1)\cap A_2$ that is outermost in $D_1$ and
$\partial$--parallel in $A_2$.  Since $\alpha'$ is outermost,
$\alpha'$ and a subarc of $\beta_1$, say $\beta_1'$, bound a subdisk
$d_1$ of $D_1$ and $d_1\cap A_2=\alpha'$.  Since $\alpha'$ is
$\partial$--parallel in $A_2$, $\alpha'$ and a subarc of $\partial
A_2$, say $\beta_2'$, bound a subdisk $d_2$ of $A_2$.  Moreover, since
$d_1\cap A_2=\alpha'$, $\beta_1'\cup\beta_2'$ bounds an embedded bigon
$\Delta'$ in $\partial M$ and $d_1\cup d_2\cup\Delta'$ bounds a
football region, which we denote by $X'$.

If $\int(d_2)\cap D_1=\emptyset$, then either $X'\subset X$ or
$\int(X)\cap \int(X')=\emptyset$.  Since the football region $X$ is
assumed to be innermost, $X'$ does not lie in $X$.  Moreover, since
$\partial A_1\cup\partial A_2$ does not form any trivial bigon, both
football regions $X$ and $X'$ must contain the vertex of the
triangulation. This means that $\int(X)\cap \int(X')\ne\emptyset$.  Thus
$\int(d_2)\cap D_1\ne\emptyset$.

Let $\alpha''\subset d_2\cap D_1$ be an outermost intersection arc in
$d_2$.  We use $e_2$ to denote the subdisk of $d_2$ bounded by
$\alpha''$ and $\beta_2'$ ($e_2\cap D_1=\alpha''$).  As
$\alpha''\subset D_1$, the arc $\alpha''$ and a subarc of $\beta_1$
bound a subdisk of $D_1$ which we denoted by $e_1$.  As before, $e_1$,
$e_2$ and a bigon in $\partial M$ bound another football region, which
we denote by $X''$.  Since $e_1\subset D_1$ and $e_2\cap
D_1=\alpha''$, if $e_2$ lies in the football region $X$, then
$X''\subset X$ contradicting the assumption the $X$ is innermost.
Similarly, if $e_2$ is outside $X$, then since $e_2\cap D_1=\alpha''$,
$X''$ must be outside $X$ and $X''\cap \int(X)=\emptyset$.  As before,
this is also impossible because by our assumptions every football
region must contain the vertex of the triangulation, which implies
$X''\cap \int(X)\ne\emptyset$.
\end{proof}

\begin{corollary}\label{Cess}
Let $\alpha$ be an arc component of $A_1\cap A_2$ and suppose $\alpha$
is $\partial$--parallel in $A_1$.  Then the following are true.
\begin{enumerate}
\item $\alpha$ must be outermost in $A_1$.
\item $\alpha$ must be an essential arc in $A_2$.
\end{enumerate}
\end{corollary}
\begin{proof}
We first prove that if $\alpha$ is outermost in $A_1$ then $\alpha$
must be an essential arc in $A_2$. Suppose otherwise that $\alpha$ is
$\partial$--parallel in $A_2$.  Since $\alpha$ is outermost in $A_1$,
the two subdisks of $A_1$ and $A_2$ cut off by $\alpha$ form an
embedded disk and bound a football region, which contradicts
\fullref{LFB}.

Since $\alpha$ is $\partial$--parallel in $A_1$, $\alpha$ and a subarc
of $\partial A_1$ bound a subdisk $D$ of $A_1$.  Suppose $\alpha$ is
not outermost.  Then we can choose $\alpha$ so that every component of
$\int(D)\cap A_2$ is outermost in $A_1$.  Let $\alpha_1,\dots,\alpha_k$
be the components of $\int(D)\cap A_2$.  Since each $\alpha_i$ is
outermost, by the argument above, every $\alpha_i$ is an essential arc
in $A_2$.  This is an immediate contradiction to \fullref{Louter}.

Part (2) follows from part (1) and the argument above.
\end{proof}

\begin{figure}[ht!]
\begin{center}
\labellist\small \pinlabel (a) [t] at 132 503 \pinlabel (b) [t] at 404
503 \pinlabel $\Gamma_1$ at 45 582 \pinlabel $\Gamma_1$ at 320 582
\pinlabel $\Gamma_2$ at 111 597 \pinlabel $\Gamma_2$ at 368 595
\pinlabel $\Gamma_2$ at 442 564 \pinlabel $Q'$ at 73 552 \pinlabel
monogon at 177 577 \endlabellist
\includegraphics[width=4in]{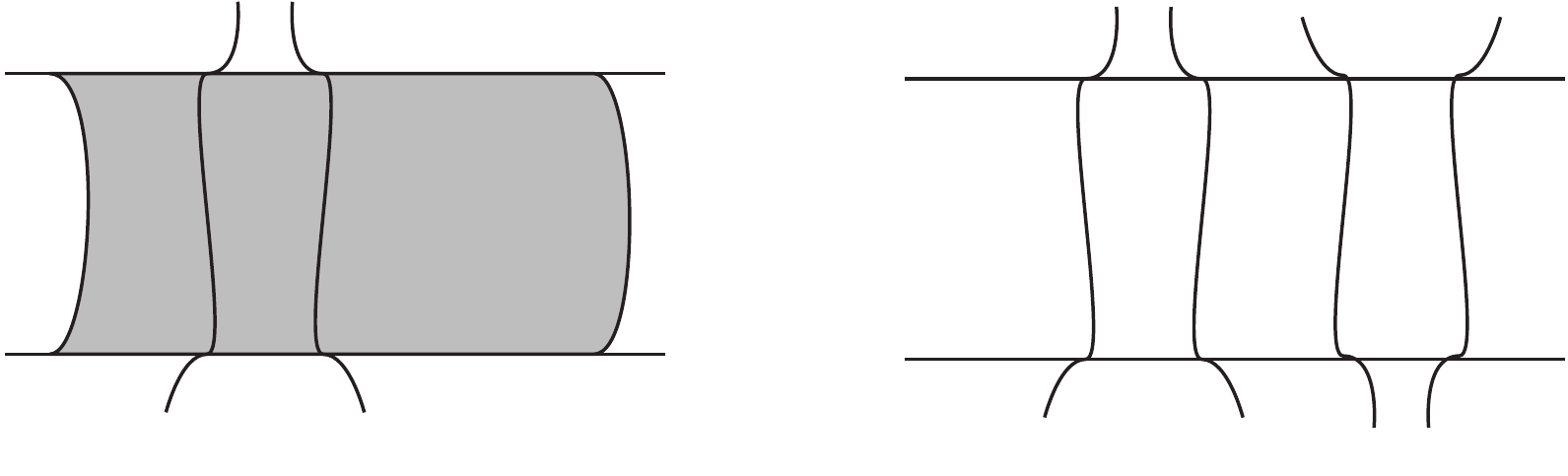}
\caption{}\label{Fmono}
\end{center}
\end{figure}

\begin{lemma}\label{Ltrivial} Suppose $A_1\cup A_2$ is bigon-efficient. 
If $A_1\cap A_2$ contains an arc that is $\partial$--parallel in $A_1$
then
\begin{enumerate}
\item every arc of $A_1\cap A_2$ is $\partial$--parallel and outermost
in $A_1$ but essential in $A_2$,
\item $A_1\cap T_2$ consists of $\partial$--compressing disks of
$A_2$,
\item every arc of $\partial A_1\cap\Gamma_2$ is of type $I$ in
$\Gamma_2$ and every arc of $\partial A_2\cap\Gamma_1$ is of type $II$
in $\Gamma_1$, see \fullref{Ftype}
\end{enumerate}
\end{lemma}
\begin{proof}
We first claim that $A_1\cap A_2$ contains no closed curve.  Suppose
otherwise $A_1\cap A_2$ contains a closed curve.  Since $A_1$ and
$A_2$ are incompressible by \fullref{L01}, every closed curve in
$A_1\cap A_2$ is either essential in both $A_1$ and $A_2$ or trivial
in both $A_1$ and $A_2$.  Since $A_1\cup A_2$ is bigon-efficient, a
closed curve in $A_1\cap A_2$ is essential in both annuli. This
implies that every arc component of $A_1\cap A_2$ is
$\partial$--parallel in both $A_1$ and $A_2$, a contradiction to
\fullref{Cess}.

Suppose $A_1\cap A_2$ contains an arc which is essential in $A_1$ and
let $\gamma_1,\dots,\gamma_k$ be all the components of $A_1\cap A_2$
that are essential in $A_1$.  Then $\gamma_1,\dots,\gamma_k$ cut $A_1$
into a collection of rectangles $R_1,\dots,R_k$ and we can suppose
$R_i$ is the rectangle between $\gamma_i$ and $\gamma_{i+1}$ (setting
$\gamma_{k+1}=\gamma_1$).  In other words, $\gamma_i$ and
$\gamma_{i+1}$ are two opposite edges of $R_i$ and the other two edges
of $R_i$ are subarcs of $\partial A_1$.

Since $A_1\cap A_2$ contains an arc trivial in $A_1$, at least one
$R_i$ contains other arcs of $A_1\cap A_2$. Let $\alpha_1,\dots,
\alpha_m$ be the components of $\int(R_i)\cap A_2$.  By our
construction of $R_i$, each $\alpha_j$ is $\partial$--parallel in
$A_1$.  By \fullref{Cess}, each $\alpha_j$ is
$\partial$--parallel and outermost in $A_1$.  Hence each $\alpha_j$
and a subarc of $\partial A_1$ bound a disk $\Delta_j$ in $R_i$ and
these $\Delta_j$'s are pairwise disjoint.  Moreover, each $\Delta_j$
is a $\partial$--compressing disk of $A_2$, in particular
$\Delta_j\subset T_2$.  This implies that $\partial A_2$ and the arcs
$\partial\Delta_j\cap\Gamma_2$ naturally deform into a Reeb train
track.  By \fullref{PReeb}, $\Gamma_2$ contains the vertex of
the triangulation.

Let $P$ and $M'$ be the closures of $R_i-\bigcup_{j=1}^m\Delta_j$ and
$M-T_2$ respectively.  So $P$ is a disk properly embedded in $M'$.
Let $P'$ be the disk bounded by $\partial P$ in $\partial M'$.  We may
consider $A_2$ as an annulus in $\partial M'$ and $P'\cap
A_2\ne\emptyset$.  Similar to the proof of \fullref{Louter},
$\partial A_2$ cuts $P'$ into a collection of disks and there are at
least two outermost such disks.  If an outermost disk $\Delta$ lies in
$\partial M'-\int(A_2)=\partial M-\int(\Gamma_2)$, as in the proof of
\fullref{Louter}, $\Delta$ must contain the vertex, which
contradicts the previous conclusion that the vertex lies in
$\Gamma_2$.  This means that every outermost disk in $P'-\partial A_2$
lies in $A_2$.  Since each $\alpha_j$ ($j=1,\dots,m$) is essential in
$A_2$, this implies that there are exactly two outermost disks and
both $\gamma_i$ and $\gamma_{i+1}$ must be $\partial$--parallel arcs
in $A_2$.

Let $\beta_i$ and $\beta_{i+1}$ be subarcs of $\partial A_2$ such that
$\partial\gamma_i=\partial\beta_i$,
$\partial\gamma_{i+1}=\partial\beta_{i+1}$, and $\gamma_i\cup\beta_i$
and $\gamma_{i+1}\cup\beta_{i+1}$ bound subdisks $\delta_i$ and
$\delta_{i+1}$ of $A_2$ respectively.  By \fullref{Cess},
$\beta_i$ and $\beta_{i+1}$ must both be outermost in $A_2$ and
$\delta_i$ and $\delta_{i+1}$ are disjoint $\partial$--compressing
disks for $A_1$.  This implies that $\beta_i$ and $\beta_{i+1}$ are of
type $I$ in $\Gamma_1$.

Note that $R_i\cup\delta_i\cup\delta_{i+1}$ is a disk properly
embedded in $M$.  Moreover, $\partial A_1\cup\beta_i\cup\beta_{i+1}$
naturally deforms into a Reeb train track and
$\partial(R_i\cup\delta_i\cup\delta_{i+1})$ deforms into a bigon in
the Reeb train track.  Let $Q'$ be the disk bounded by
$\partial(R_i\cup\delta_i\cup\delta_{i+1})$ in $\partial M$, see the
shaded region in \fullref{Fmono}(a) for a picture of $Q'$. Clearly
$Q'\subset\Gamma_1$.  As above, we say a disk in $Q'-\int(\Gamma_2)$ is
outermost if its boundary consists of an arc from $\partial A_1$ and
an arc from $\partial A_2$.  As in the proof of \fullref{Louter},
any outermost disk must contain the vertex of the triangulation. Since
$\Gamma_2$ contains the vertex, $Q'-\int(\Gamma_2)$ contains no
outermost disk.  This implies that $Q'\cap\Gamma_2$ consists of
rectangles which naturally deform into bigons in the Reeb annulus
$\Gamma_2$, see \fullref{Fmono}(a) for a picture.  As shown in
\fullref{Fmono}(a), at least one component of $Q'-\int(\Gamma_2)$ is
a monogon (after deforming into a train track).  Since a monogon
contains the vertex of the triangulation, this implies that the vertex
of the triangulation lies outside $\Gamma_2$, a contradiction.  This
proves part (1).

Part (2) is an immediate corollary of part (1).

Part (1) also implies that every arc of $\Gamma_2\cap\partial A_1$ is
of type $I$ in $\Gamma_2$ and $\partial A_2\cup(\Gamma_2\cap\partial
A_1)$ forms a standard Reeb train track.  Now we consider $\partial
A_2\cap\Gamma_1$.

As above, since $\Gamma_2$ contains the vertex, $\Gamma_1-\Gamma_2$
has no outermost disk (an outermost disk is a component with a
boundary edge in $\partial A_1$ and a boundary edge in $\partial
A_2$). This implies that every arc in $\partial A_2\cap\Gamma_1$ is an
essential arc in $\Gamma_1$.  Since $\partial
A_2\cup(\Gamma_2\cap\partial A_1)$ form a standard Reeb train track,
as shown in \fullref{Fmono}(b), every arc of $\partial
A_2\cap\Gamma_1$ must be of type $II$ in $\Gamma_1$.
\end{proof}

\begin{lemma}\label{Ltype}
Suppose $A_1\cap A_2$ is bigon-efficient and $A_1\cap
A_2\ne\emptyset$.  Then no arc component of $A_1\cap A_2$ is essential
in both $A_1$ and $A_2$.
\end{lemma}
\begin{proof}
Suppose there is an arc component of $A_1\cap A_2$ that is essential
in both $A_1$ and $A_2$.  As in the proof of \fullref{Ltrivial},
$A_1\cap A_2$ contains no closed curve.  If there is a component of
$A_1\cap A_2$ that is trivial in $A_1$ then by \fullref{Ltrivial}
every component of $A_1\cap A_2$ is trivial in $A_1$.  Thus every arc
of $A_1\cap A_2$ must be essential in both $A_1$ and $A_2$.

So $A_1\cap A_2$ cuts both $A_1$ and $A_2$ into a collections of
rectangles.  Let $R$ be a component of $A_1\cap T_2$.  Two opposite
boundary edges of the rectangle $R$ are essential arcs in $A_2$ and
the other two edges of $\partial R$, denoted by $\gamma_1$ and
$\gamma_2$, are properly embedded in $\Gamma_2$.  Since $R$ is a disk
properly embedded in the solid torus $T_2$, both $\gamma_1$ and
$\gamma_2$ must be $\partial$--parallel in $\Gamma_2$.  Moreover,
since each arc in $A_1\cap A_2$ is essential in both $A_1$ and $A_2$,
$\partial\gamma_1$ and $\partial\gamma_2$ lie in different components
of $\partial\Gamma_2$.  Thus $\gamma_1$ and $\gamma_2$ and two subarcs
of $\partial A_2$ (from different components of $\partial A_2$) bound
two disjoint disks $d_1$ and $d_2$ in $\partial M$ respectively.
After naturally deforming $\partial A_1\cup\partial A_2$ into a train
track, $d_1$ and $d_2$ become bigons or monogons.  Since $\partial
A_1\cap\partial A_2$ is bigon-efficient, every bigon contains the
vertex of the triangulation.  Since every monogon also must contain
the vertex, this contradicts that $d_1$ and $d_2$ are disjoint and
there is only one vertex in the triangulation.
\end{proof}

\begin{corollary}\label{CtypeI}
Suppose $A_1\cap A_2$ is bigon-efficient. Suppose $A_1$ is of type $I$
relative to $A_2$. Then every arc of $\Gamma_1\cap\partial A_2$ is of
type $I$ in $\Gamma_1$ and every arc of $\Gamma_2\cap\partial A_1$ is
of type $II$ in $\Gamma_2$.\end{corollary}
\begin{proof} By \fullref{Ltype}, no arc component of $A_1\cap A_2$ is essential in both $A_1$ and $A_2$.  Now the corollary follows from \fullref{Ltrivial}.
\end{proof}

Next we study the intersection patterns of 3 normal annuli carried by
$B$.

\begin{lemma}\label{LA3}
Let $A_1$, $A_2$ and $A_3$ be pairwise bigon-efficient normal annuli
carried by $B$.  Suppose $A_1$ is of type $I$ relative to $A_2$ and
$\partial A_1\cap\partial A_3\ne\emptyset$.  Then,
\begin{enumerate}
\item $A_1$ must be of type $I$ relative to $A_3$ and
\item $\partial A_2\cap\partial A_3=\emptyset$.
\end{enumerate}
\end{lemma}
\begin{proof}
Since $\partial A_1\cap\partial A_3\ne\emptyset$, by
Lemmas~\ref{Ltype} and \ref{Ltrivial}, either $A_1$ is of type $I$
relative to $A_3$ or $A_3$ is of type $I$ relative to $A_1$.  Suppose
part (1) is not true and $A_3$ is of type $I$ relative to $A_1$.  So
by \fullref{Ltrivial} and \fullref{PReeb}, $\Gamma_3$
contains the vertex of the triangulation.  Moreover, since $A_1$ is of
type $I$ relative to $A_2$, both $\Gamma_1$ and $\Gamma_3$ contain the
vertex.

Let $R$ be the component of $\Gamma_1\cap\Gamma_3$ that contains the
vertex of the triangulation.  By \fullref{Ltrivial}, $\partial
A_1\cap\Gamma_3$ consists of type $I$ arcs in $\Gamma_3$, so $R$ is a
quadrilateral that naturally deforms into a bigon.  Two opposite edges
of $\partial R$, denoted by $r_1$ and $r_2$, are components of
$\Gamma_1\cap\partial A_3$.  By \fullref{Ltrivial}, $r_1$ and $r_2$
are type $II$ arcs in $\Gamma_1$, see \fullref{F3}(a) for a picture
of $R$.  Let $r_3$ and $r_4$ be the other two edges of $R$.  Hence,
$r_3\cup r_4$ are two components of $\partial A_1\cap\Gamma_3$ and
$r_3$ and $r_4$ are of type $I$ in $\Gamma_3$.

Since $A_1$ is of type $I$ relative to $A_2$, every component of
$\partial A_2\cap\Gamma_1$ is of type $I$ in $\Gamma_1$ and $\partial
A_1\cup (\partial A_2\cap\Gamma_1)$ forms a standard Reeb train track.

\medskip

\noindent\textbf{Case 1}\qua  $(\partial A_2\cap\Gamma_1)\bigcap (r_1\cup
r_2)=\emptyset$

\medskip

If a component of $\partial A_2\cap\Gamma_1$ lies outside $R$, as
shown in \fullref{F3}(a), it creates a monogon region outside $R$.
Since any monogon region contains the vertex, this contradicts that
$R$ contains the vertex.  Thus $\partial A_2\cap\Gamma_1\subset R$.

\begin{figure}[ht!]
\begin{center}
\labellist \small \pinlabel* (a) [t] at 170 559 \pinlabel* (b) [t] at
367 559 \pinlabel* (c) [t] at 170 423 \pinlabel* (d) [t] at 367 423
\pinlabel* (e) [t] at 282 317 \pinlabel $\Gamma_1$ at 112 629
\pinlabel $\Gamma_3$ at 154 574 \pinlabel $\Gamma_3$ at 336 574
\pinlabel $\Gamma_3$ at 109 470 \pinlabel $\Gamma_3$ at 313 470
\pinlabel $\Gamma_3$ at 215 365 \pinlabel $R$ at 168 621 \pinlabel
$R'$ at 184 484 \pinlabel $R'$ at 387 489 \pinlabel $R'$ at 285 365
\hair 2pt \pinlabel $r_1$ [r] at 143 622 \pinlabel $r_2$ [r] at 194
621 \pinlabel $r_3$ [b] at 170 647 \pinlabel $r_3$ [b] at 354 647
\pinlabel $r_3$ [l] at 146 485 \pinlabel $r_3$ [l] at 344 490
\pinlabel $r_3$ [r] at 250 386 \pinlabel* $r_4$ [r] at 206 483
\pinlabel* $r_4$ [r] at 405 487 \pinlabel $r_4$ [r] at 312 346
\pinlabel $r_4$ [t] at 174 593 \pinlabel $r_4$ [t] at 357 593
\pinlabel $\alpha$ [t] at 172 474 \pinlabel $\alpha$ [t] at 362 462
\pinlabel $X$ at 334 620 \pinlabel $Y$ at 378 620 \pinlabel monogon
<1pt,0pt> at 219 621 \pinlabel monogon <1pt,0pt> at 173 451 \pinlabel
monogon <2pt,-2pt> at 369 473 \endlabellist
\includegraphics[width=4in]{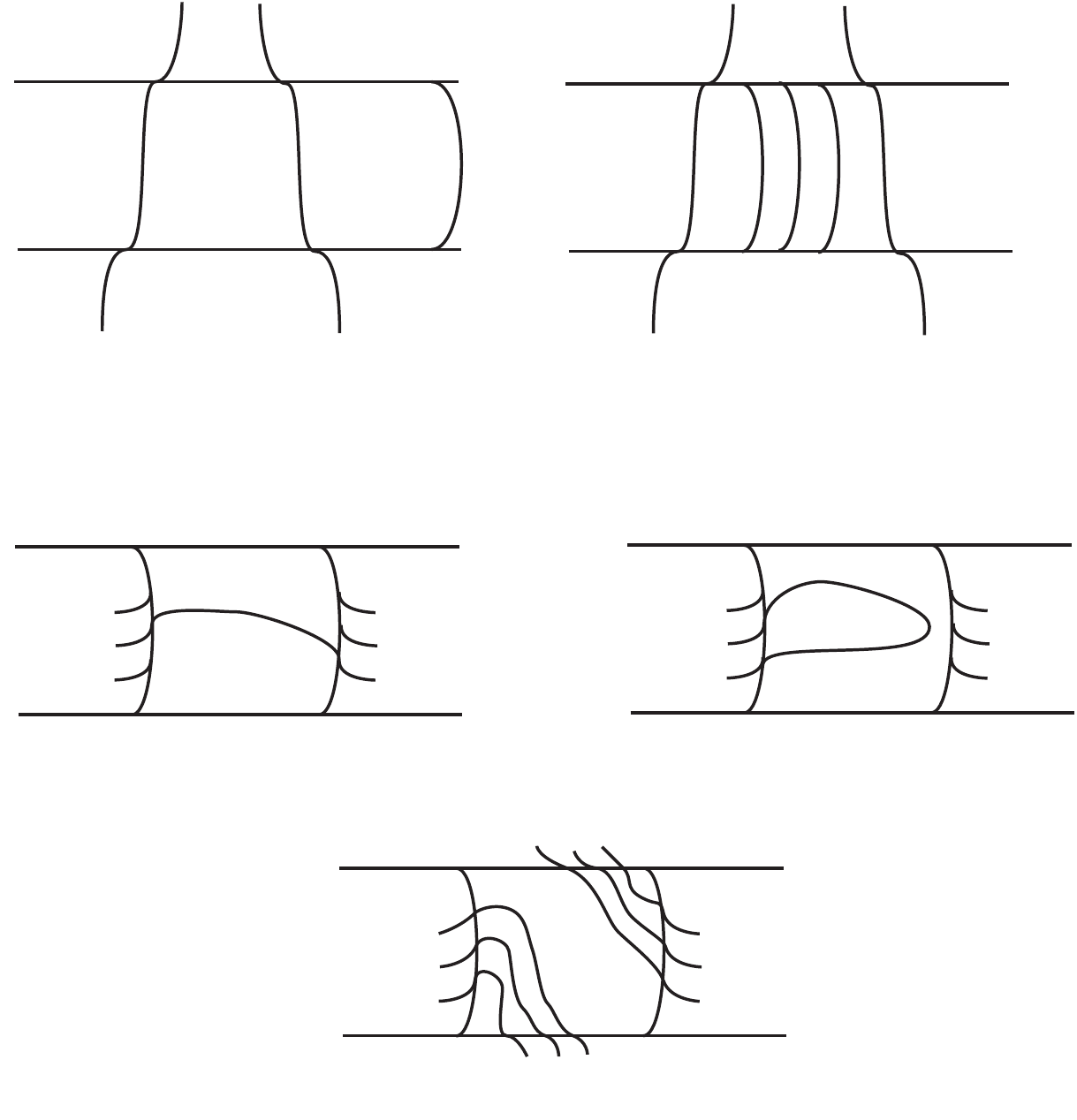}
\caption{}\label{F3}
\end{center}
\end{figure}

Next we view $R$ as a quadrilateral in $\Gamma_3$. Hence $r_3$ and
$r_4$ are type $I$ arcs in $\Gamma_3$.  Each component of $\partial
A_2\cap R$ is an arc with one endpoint in $r_3$ and the other endpoint
in $r_4$.  Moreover, as shown in \fullref{F3}(b), after deforming
into a train track, $\partial A_2\cap R$ cuts $R$ into a monogon
region $X$, a 3--prong triangle $Y$, and a collection of bigons.

Now we consider the disk $R'=\overline{\Gamma_3-R}$.  We first
consider the possibility that there is an arc $\alpha$ of $\partial
A_2\cap R'$ with both endpoints in $r_3\cup r_4$.  Note that $\partial
A_2\cap R$ is as shown in \fullref{F3}(b), so this configuration
fixes the switching direction of $\partial\alpha$ in the train track.
There are two cases to consider: (1) both endpoints of $\alpha$ lie in
$r_3$ (or $r_4$) and (2) one endpoint of $\alpha$ lies in $r_3$ and
the other lies in $r_4$.  As shown in \fullref{F3}(c) and (d), in
either case, $\alpha$ produces a monogon in $R'$, which means the
vertex of the triangulation lies in $R'$ and contradicts the
assumption that $R$ contains the vertex.  Thus, every component of
$\partial A_2\cap R'$ has one endpoint in $r_3\cup r_4$ and the other
endpoint in $\partial\Gamma_3\cap\partial R'$.

After deforming into a train track, $R'$ becomes a bigon. Since $R'$
does not contain the vertex, $\partial A_2$ cuts $R'$ into a
collection of disks, each of which becomes a bigon after deformed into
a train track.  Because of the switching direction of the train track
at $\partial A_2\cap(r_3\cup r_4)$, as shown in \fullref{F3}(e),
the arcs with an endpoint in $r_4$ must have the same configuration.
Otherwise, these arcs would create a monogon in $R'$.  Furthermore,
since every arc of $\partial A_2\cap\Gamma_3$ is essential in
$\Gamma_3$ by \fullref{Ltrivial}, the arcs with an endpoint in $r_3$
must also have the same configuration as shown in \fullref{F3}(e).
In other words, \fullref{F3}(e) is the only possible configuration
for $\partial A_2\cap R'$.

As shown in \fullref{Fspiral}(a) and (b), given a component
$\alpha$ of $\partial A_3$ and any arc $\beta$ intersecting $\alpha$,
there are essentially two different switching directions at
$\alpha\cap\beta$ along $\alpha$.  By examining the switching
directions of the train track at $\partial A_2\cap\partial A_3$ in
$\partial R'$ along $\partial A_3$ as shown in \fullref{F3}(e), we
can see that each component of $\partial A_2\cap\Gamma_3$ must be of
type $II$ in $\Gamma_3$.  Moreover, as shown in \fullref{F3}(e),
the argument above implies that the switching directions (of the train
track) at the intersection points of $\partial A_2$ with any component
of $\partial A_3$ are all the same.  However, by part (3) of
\fullref{Ltrivial}, the conclusion that $\partial A_2\cap\Gamma_3$
contains a type $II$ arc in $\Gamma_3$ implies that $\partial
A_3\cap\Gamma_2$ consists of type $I$ arcs in $\Gamma_2$.  This means
that there are two arcs of $\partial A_2\cap\Gamma_3$, similar to the
$r_1$ and $r_2$ in \fullref{F3}(a), whose endpoints on a component
of $\partial\Gamma_3$ have opposite switching direction.  This
contradicts the previous conclusion (as depicted in
\fullref{F3}(e)) that all the switching directions at such points
are the same.

\medskip

\noindent\textbf{Case 2}\qua  $(\partial A_2\cap\Gamma_1)\bigcap (r_1\cup
r_2)\ne\emptyset$

\medskip

We will perform some normal isotopies so that $(\partial
A_2\cap\Gamma_1)\bigcap (r_1\cup r_2)=\emptyset$ after the isotopies.

Let $\alpha_i\subset\partial A_i$ ($i=1,2,3$) be 3 arcs intersecting
each other and forming a triangle $\Delta$ as shown in
\fullref{Fbigon}.  Suppose $\Delta$ naturally deforms into a bigon
and $\Delta$ does not contain the vertex of the triangulation.  Then,
as shown in \fullref{Fbigon}, the isotopy on $\alpha_3$, fixing
$\alpha_1$ and $\alpha_2$, is a normal isotopy.  Next we will fix
$\partial A_1\cup\partial A_3$ and perform some isotopies as in
\fullref{Fbigon} so that $(\partial A_2\cap\Gamma_1)\bigcap
(r_1\cup r_2)=\emptyset$ after the isotopies.  Each isotopy pushes an
intersection point of $\partial A_2\cap\partial A_3$ out of
$\Gamma_1$.

\begin{figure}[ht!]
\begin{center}
\labellist\small \pinlabel {deform into train track} [rb] at 212 529
\pinlabel {normal isotopy} [lb] at 297 529 \pinlabel {$\Delta$} at 255
587 \endlabellist \includegraphics[width=4in]{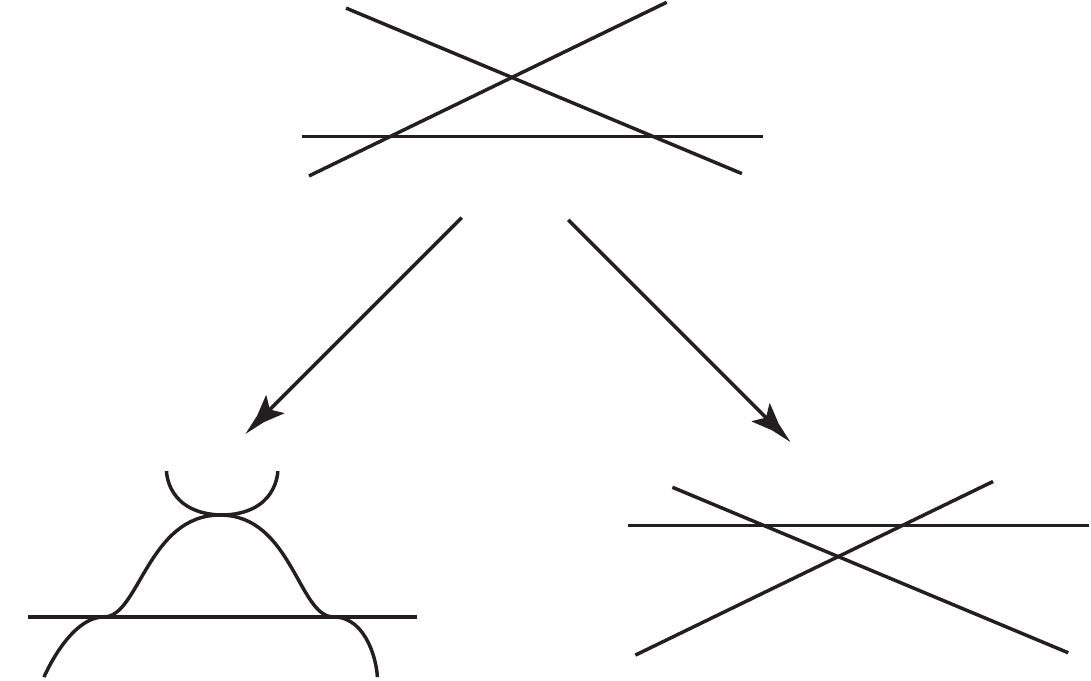}
\caption{}\label{Fbigon}
\end{center}
\end{figure}

Let $\alpha$ be a component of $\partial A_2\cap\Gamma_1$ and suppose
$\alpha\cap (r_1\cup r_2)\ne\emptyset$. Let $\alpha_1$ and $\alpha_2$
be the closure of the components of $\alpha-(r_1\cup r_2)$ that
contain $\partial\alpha$.  So $\alpha_i$ ($i=1,2$) has one endpoint in
$\partial A_1$ and the other endpoint in $r_1\cup r_2$.  Thus
$\alpha_1$ and $\alpha_2$ are the edges of two triangles $\Delta_1$
and $\Delta_2$ respectively formed by $\partial A_1$, $\partial A_2$
and $\partial A_3$. Since the two endpoints of $\alpha$ lie in
different components of $\partial\Gamma_1$, $\Delta_1$ and $\Delta_2$
are not nested.  Without loss of generality, we may assume each
$\Delta_i$ is innermost, ie, $\int(\Delta_i)\cap (\partial
A_1\cup\partial A_2\cup\partial A_3)=\emptyset$ for both $i=1,2$.
After deforming $\partial A_1\cup\partial A_2\cup\partial A_3$ into a
train track, each $\Delta_i$ becomes either a bigon or a monogon.
Since each monogon contains the vertex of the triangulation, at least
one $\Delta_i$ is a bigon that does not contain the vertex.  Hence a
normal isotopy on $\partial A_2$, as shown in \fullref{Fbigon}
pushes an intersection point of $\partial A_2\cap (r_1\cup r_2)$ out
of $\Gamma_1$.  So after finitely many such normal isotopies,
$(\partial A_2\cap\Gamma_1)\bigcap (r_1\cup r_2)=\emptyset$ and we can
apply Case 1 to obtain a contradiction.

Therefore, $A_1$ is of type $I$ relative to both $A_2$ and $A_3$ and
part (1) of the lemma holds.  If $\partial A_2\cap\partial
A_3\ne\emptyset$, then by \fullref{Ltrivial} and \fullref{Ltype},
either $A_2$ is of type $I$ relative to $A_3$, or $A_3$ of type $I$
relative to $A_2$.  Both possibilities contradict part (1), since
$A_1$ is of type $I$ relative to both $A_2$ and $A_3$.
\end{proof}

\section{Boundary curves}\label{Scurve}

Suppose $A_1\cap A_2$ is bigon-efficient.  If $A_1\cap A_2$ contains a
closed curve, by \fullref{Ltrivial}, all the components of $A_1\cap
A_2$ must be closed essential curves.  After performing canonical
cutting and pasting along these curves, we get a pair of disjoint
annuli $A_1'$, $A_2'$ and a possible collection of tori $T$.  Clearly,
$A_1+A_2=A_1'+A_2'+T$.  In particular, $weight(A_1'+A_2')\le
weight(A_1+A_2)$.

Let $A_1,\dots, A_n$ be a fixed set of normal annuli carried by $B$.
We consider $m_i$ parallel copies of $A_i$ ($i=1,\dots,n$).  Then we
can perform the isotopy and cutting and pasting above on each pair of
the $\sum_{i=1}^nm_i$ annuli, so that each pair of resulting set of
annuli are bigon-efficient and have no closed intersection curve.  So
there is a set of normal annuli $\mathcal{A}$ such that for any set of
nonnegative integers $m_i$, there is a collection of annuli
$A_1',\dots,A_k'$ in $\mathcal{A}$ such that
\begin{enumerate}
\item $\sum_{i=1}^n m_iA_i$=$T+\sum_{i=1}^k m_i'A_i'$, where $T$ is a
collection of normal tori.
\item $A_1',\dots,A_k'$ are pairwise bigon-efficient
\item $A_i'\cap A_j'$ contains no closed curve for any $i\ne j$.
\end{enumerate}
We claim that one can choose $\mathcal{A}$ to be a finite set of
annuli.  Let $\partial\mathcal{A}$ be the set of boundary curves of
all possible normal annuli resulting from the normal isotopies and
canonical cutting and pasting as above (among all possible $m_i$'s).
Since the operations that make $A_k\cap A_j$ bigon-efficient, when
restricted to $\partial M$, are simply cutting and pasting on bigons
in $\partial M$, $\partial\mathcal{A}$ is a finite set of normal
curves. Now suppose there is an infinite set of normal annuli, denoted
by $\mathcal{D}$, in $\mathcal{A}$ with the same pair of boundary
curves.  Then by the normal surface theory, there must be two annuli
in $\mathcal{D}$, say $A_i'$ and $A_j'$ such that $A_j'=T'+A_i'$ where
$T'$ is a collection of normal tori.  This means that $A_j'$ is
redundant as we can use $T'+A_i'$ instead.  Therefore, we may choose
$\mathcal{A}$ to be a finite set and there is an algorithm to find all
the annuli in $\mathcal{A}$ using normal surface theory.

We are mainly interested in the boundary curves.  In the conclusion
(1) above, clearly $\sum_{i=1}^n m_i\partial A_i$=$\sum_{i=1}^k
m_i'\partial A_i'$ in $\partial M$.

Since each $A_i'$ is also a normal annulus, $A_i'$ is
$\partial$--parallel in $\partial M$.  We use $\Gamma_i'$ to denote
the subannulus of $\partial M$ isotopic to $A_i'$ and with
$\partial\Gamma_i'=\partial A_i'$.

Let $\mathcal{S}$ be the union of a fixed set of pairwise disjoint
compact surfaces carried by $B$. Since $\partial\mathcal{S}$ is
carried by the train track $\partial B$ and $\partial B$ does not
carry any trivial circle, every component of $\partial\mathcal{S}$ is
an essential normal curve.  Next we consider $\partial\mathcal{S}+\sum
m_i\partial A_i'$.  Our goal is to prove the following lemma.

\begin{lemma}\label{Lcurve}
Let $\mathcal{S}$ and $A_i'$ be as above.  Then the diameter of the
set $\{\partial\mathcal{S}+\sum m_i\partial A_i'\}$ (for all
nonnegative integers $m_i$) in the curve complex $\mathcal{C}(F)$ is
bounded.
\end{lemma}

Suppose $A_1'$ and $A_2'$ are of type $I$ relative to $A_i'$ and
$A_j'$ respectively.  If $\partial A_1'\cap\partial A_2'\ne\emptyset$,
then by \fullref{Ltrivial} and \fullref{Ltype}, one of $A_1'$ and
$A_2'$ is of type $II$ relative the other, contradicting
\fullref{LA3}.  Thus $\partial A_1'\cap\partial A_2'=\emptyset$.
Since both $\Gamma_1'$ and $\Gamma_2'$ contain the vertex of the
triangulation by \fullref{PReeb} and since $A_1'\cap A_2'$
contains no closed curve, $\Gamma_1'$ and $\Gamma_2'$ must be nested
and $\partial A_1'$ must be normally isotopic to $\partial A_2'$.
Thus we have $m_1\partial A_1'+m_2\partial A_2'=(m_1+m_2)\partial
A_1'$.

We say $A_i'$ is of type $I$ if $A_i'$ is of type $I$ relative to one
of $A_1',\dots,A_k'$.  The argument above implies that the boundary of
all the type $I$ annuli are normally parallel.  Moreover, by
\fullref{Ltrivial} and \fullref{LA3}, those annuli among
$A_1',\dots,A_k'$ that are not of type $I$ are pairwise disjoint.

Next we will only focus on the boundary curves of
$A_1',\dots,A_k'$. If no $A_i'$ is of type $I$, then
\fullref{Ltrivial} and \fullref{Ltype} imply that these $A_i'$'s
are mutually disjoint.  Suppose $A_1'$ a type $I$ annulus.  Since the
boundary of other type $I$ annuli are normally parallel to $\partial
A_1'$, without loss of generality, we may assume $A_1'$ is the only
type $I$ annulus in $A_1',\dots,A_k'$.  By \fullref{Ltrivial} and
\fullref{LA3}, this implies $A_2'\dots,A_k'$ are pairwise
disjoint. Let $\gamma_i$ be a component of $\partial A_i$ and $k_i$
the number of intersection points of $\gamma_i$ with
$\partial\mathcal{S}$.

\begin{lemma}\label{LH}
The distance between $\gamma_j$ ($j\ne 1$) and
$\partial\mathcal{S}+\sum_{i=2}^k m_i\partial A_i'$ is at most
$2+2\log_2 k_j$.
\end{lemma}
\begin{proof}
By our earlier assumptions, $A_2',\dots, A_k'$ are mutually disjoint.
So $\sum_{i=2}^k m_i\partial A_i'$ is a union of disjoint curves and
we may regard $\gamma_j$ ($j\ne 1$) as a component of $\sum_{i=2}^k
m_i\partial A_i'$.  Since the number of intersection points of
$\gamma_j$ with $\partial\mathcal{S}$ is $k_j$, the intersection
number of $\gamma_j$ and $\partial\mathcal{S}+\sum_{i=2}^k m_i\partial
A_i'$ is at most $k_j$.  Now it is clear that \fullref{LH} follows
from \cite[Lemma 2.1]{H}, which says that the distance between any two
curves with intersection number $k$ is at most $2+2\log_2 k$.
\end{proof}

Note that \fullref{LH} implies \fullref{Lcurve} in the case that
no $A_i'$ is of type $I$.

\begin{lemma}\label{LH2}
If there is some $\partial A_j'$ ($j\ne 1$) disjoint from $\partial
A_1'$, then the distance between $\gamma_j$ and
$\partial\mathcal{S}+\sum_{i=1}^k m_i\partial A_i'$ is at most
$2+2\log_2 k_j$.
\end{lemma}
\begin{proof}
The proof is identical to that of \fullref{LH}.  Since $A_2',\dots,
A_k'$ are mutually disjoint, $\gamma_j$ can be viewed as a component
of $\sum_{i=1}^k m_i\partial A_i'$.  So the intersection number of
$\gamma_j$ and $\partial\mathcal{S}+\sum_{i=1}^k m_i\partial A_i'$ is
at most $k_j$ and \fullref{LH2} follows from Lemma 2.1 of \cite{H}.
\end{proof}

So to prove \fullref{Lcurve}, we may assume $\partial
A_i'\cap\partial A_1'\ne\emptyset$ for each $i\ne 1$.  As $A_1'$ is of
type $I$, every component of $\partial A_i'\cap\Gamma_1'$ is a type
$I$ arc in $\Gamma_1'$.

\begin{figure}[ht!]
\begin{center}
\labellist\small \pinlabel (a) [t] at 124 579 \pinlabel (b) [t] at 299
579 \pinlabel (c) [t] at 124 368 \pinlabel (d) [t] at 299 368
\pinlabel $\beta$ [b] at 137 634 \pinlabel $\beta$ [b] at 286 634
\pinlabel {cut and paste} at 210 455 \pinlabel $\alpha_i$ [b] at 120
595 \pinlabel $\alpha_i$ [b] at 300 595 \endlabellist
\includegraphics[width=3.5in]{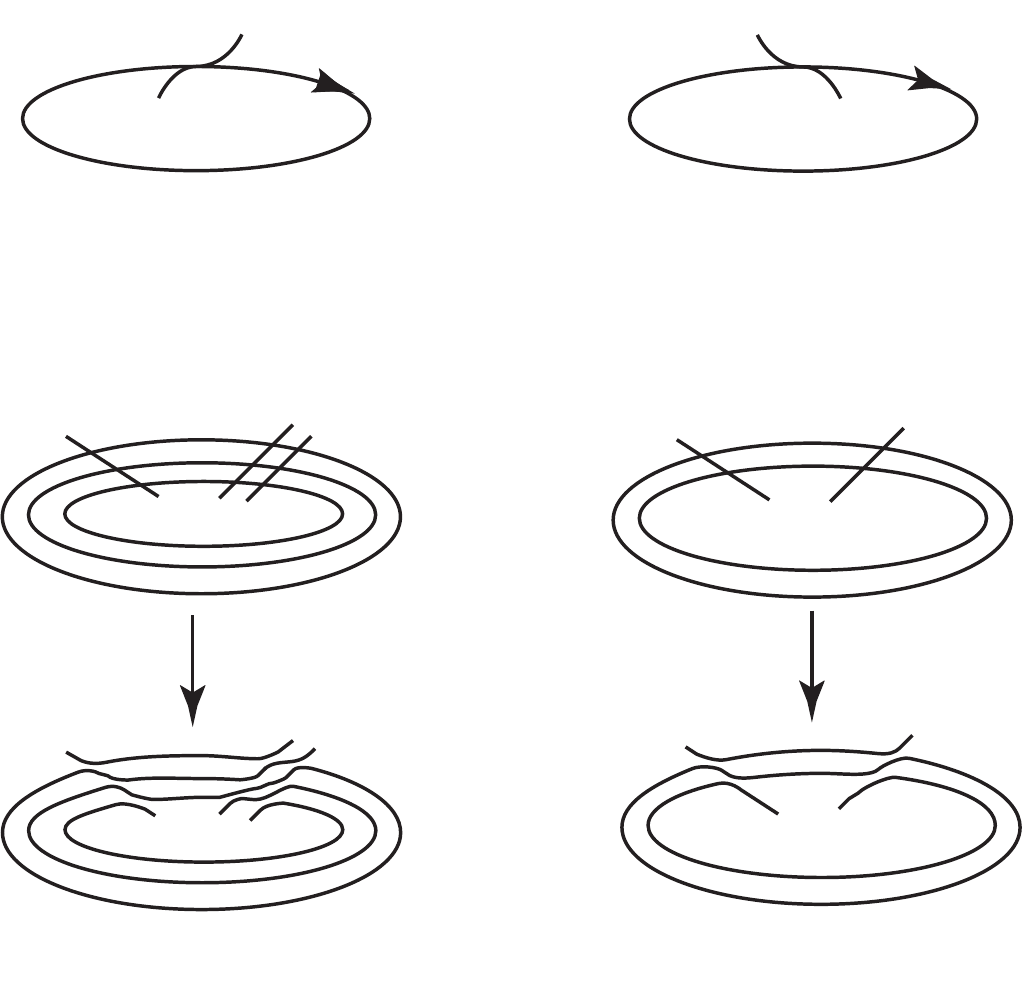}
\caption{}\label{Fspiral}
\end{center}
\end{figure}

Let $\alpha_1$ and $\alpha_2$ be the two components of $\partial
A_1'$. We fix a direction for the circle $\alpha_1$ and assign the
same direction to $\alpha_2$.  Let $\beta$ be an arc carried by
$\partial B$ and intersecting $\alpha_i$ ($i=1$ or 2) in one point. We
say $\beta$ and the point $\beta\cap\alpha_i$ are of positive
(resp.~negative) type if $\alpha_i\cup\beta$ deforms into a train
track as in \fullref{Fspiral}(a) (resp.~\fullref{Fspiral}(b)).
Note that a curve carried by the train track \fullref{Fspiral}(a)
or (b) is a spiral around $\alpha_i$.  We call a spiral carried by the
train track in \fullref{Fspiral}(a) (resp.~\fullref{Fspiral}(b))
a positive (resp.~negative) spiral.

Let $S$ be any compact surface carried by $B$ and suppose $A_1'\cap S$
contains an arc component $\gamma$.  Then there are two cases (1)
$\gamma$ is $\partial$--parallel in $A_1'$ and (2) $\gamma$ is an
essential arc in $A_1'$.  Since both $S$ and $A_1'$ are carried by the
same branched surface, as in \cite{Ha}, in either case, one endpoint
of $\gamma$ is of positive type and the other endpoint is of negative
type.  Let $P_i$ ($i=1,2$) be the number of points in $\partial
S\cap\alpha_i$ of positive type and $N_i$ the number of points in
$\partial S\cap\alpha_i$ of negative type.  The argument above implies
that $P_1+P_2=N_1+N_2$.

Let $N(\alpha_i)$ ($i=1,2$) be a small annular neighborhood of
$\alpha_i$ in $\partial M$.  We consider $\partial S+m\alpha_i$
restricted to $N(\alpha_i)$.  As depicted in \fullref{Fspiral}(c),
if $P_i\ne N_i$ and $m\ge \min\{N_i,P_i\}$, then $\partial
S+m\alpha_i$ restricted to $N(\alpha_i)$ consists of $|P_i-N_i|$
spirals and $2\min\{N_i,P_i\}$ $\partial$--parallel arcs in
$N(\alpha_i)$.  As shown in \fullref{Fspiral}(d), if $N_i=P_i$ and
$m > \min\{N_i,P_i\}$, at least one component of $\partial
S+m\alpha_i$ is parallel to $\alpha_i$ and hence we may view the
distance $d_{\mathcal{C}(F)}(\partial S+m\alpha_i, \alpha_i)\le 1$.
Without loss of generality, we may assume $P_1>N_1$.  Since
$P_1+P_2=N_1+N_2$, $P_2<N_2$. So if $m\ge \max\{N_1, P_2\}$, $\partial
S+m\partial A_1'$ has $r=P_1-N_1=N_2-P_2$ positive spirals in
$N(\alpha_1)$ and $r$ negative spirals in $N(\alpha_2)$.

\begin{figure}[ht!]
\begin{center}
\labellist\small \pinlabel (a) at 175 313 \pinlabel (b) at 444 313
\endlabellist \includegraphics[width=4in]{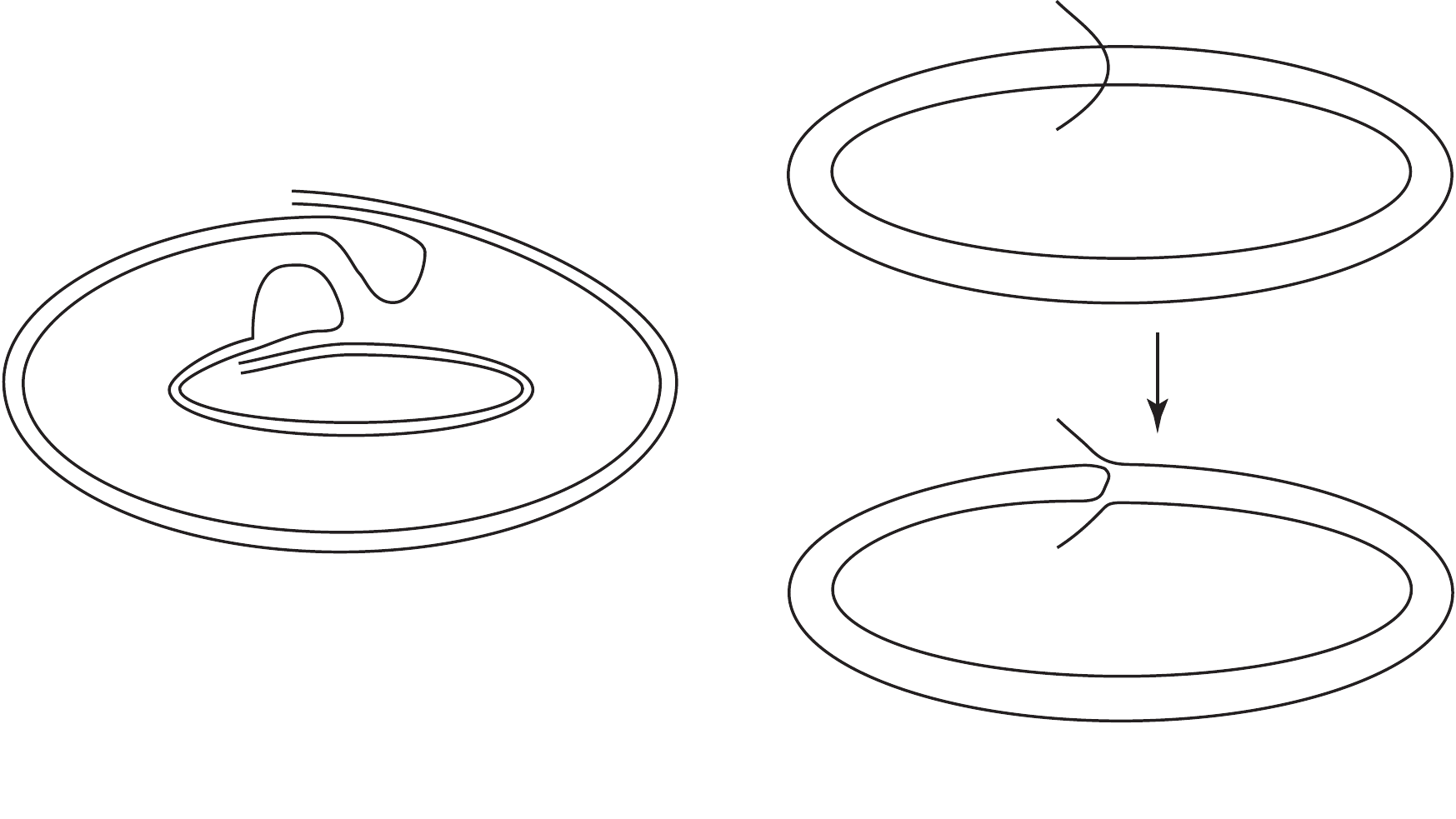}
\caption{}\label{Ftwist}
\end{center}
\end{figure}

Now we assume $m\ge \max\{N_1, P_2\}$ and consider $\partial
S+m\partial A_1'$ restricted to $N(\Gamma_1')$, which is a small
neighborhood of $\Gamma_1'$ in $\partial M$.  The positive and
negative spirals in $N(\alpha_1)$ and $N(\alpha_2)$ are connected by
some arcs of $\partial S\cap\Gamma_1'$.  First suppose two positive
spirals are connected by an arc in $\partial S\cap\Gamma_1'$. Then
this arc and the two spirals in $N(\alpha_1)$ form a monogon whose
``tail" spirals along $\alpha_1$.  Moreover, since the number of
negative spirals equals the number of positive spirals, there must be
an arc of $\partial S\cap\Gamma_1'$ connecting two negative spirals in
$N(\alpha_2)$ and hence forming another monogon, as shown in
\fullref{Ftwist}(a).  Since each monogon must contain the vertex of
the triangulation, this is a contradiction.  Thus every positive
spiral in $N(\alpha_1)$ is connected to a negative spiral in
$N(\alpha_2)$ by an arc in $\partial S\cap\Gamma_1'$.  The standard
picture of these arcs are type $I$ arcs whose two ends spiraling
around $\partial A_1'$.  Therefore, as shown in
\fullref{Ftwist}(b), $(\partial S+m\partial A_1')+\partial A_1'$ is
isotopic to $\partial S+m\partial A_1'$.

Now we assume the surface $S$ in the argument above is the resulting
surface of $\mathcal{S}+\sum_{i=2}^n m_iA_i'$ and let
$\sigma=\sum_{i=2}^n m_i$.  Clearly, there is a number $K$ depending
on $\mathcal{S}\cap A_1'$ and $A_i'\cap A_1'$, such that $K\sigma\ge
\max\{P_1,N_1, P_2, N_2\}$.  Thus by the discussion above, if $P_1\ne
N_1$ and $m_1\ge K\sigma$, the set of curves $\{\partial S+m_1\partial
A_1'\}$ are all isotopic.  Moreover, by \fullref{LH}, the set of
curves $\{\partial\mathcal{S}+\sum_{i=2}^n m_i\partial A_i'\}$ for all
$m_i$ ($i=2,\dots, n$) has bounded diameter.  As
$S=\mathcal{S}+\sum_{i=2}^n m_iA_i'$, by the argument above, if
$N_1=P_1$ \fullref{Lcurve} holds, and if $P_1\ne N_1$,
\fullref{Lcurve} holds under the condition that $m_1\ge K\sigma$.

Next we consider the case that $m_1< K\sigma$. By our assumptions,
$A_1'$ is the only type $I$ annulus and $\partial A_i'\cap\Gamma_1'$
consists of type $I$ arcs in $\Gamma_1'$. So, as in
\fullref{Ftwist}(b), $\sum_{i=1}^n m_i\partial A_i'$ is isotopic to
$\sum_{i=2}^n m_i\partial A_i'$.  Thus $\sum_{i=1}^n m_i\partial A_i'$
consists of $2\sigma$ closed curves.

Let $\omega$ be the maximal weight of $\partial A_i'$ among all $i$.
So if $m_1< K\sigma$, the total weight of $\sum_{i=1}^n m_i\partial
A_i'$ is less than $K\sigma\omega+\sigma\omega=(K+1)\sigma\omega$.
Since $\sum_{i=1}^n m_i\partial A_i'$ consists of $2\sigma$ closed
curves, there is a component $\gamma$ of $\sum_{i=1}^n m_i\partial
A_i'$ with weight less than $(K+1)\omega/2$.  Up to normal isotopy,
there are only finitely many curves with weight under $(K+1)\omega/2$.
So there is a number $K'$ such that $|\partial
\mathcal{S}\cap\gamma|\le K'$.  As in the proof of \fullref{LH}, by
a theorem of Hempel \cite{H}, the distance between $\partial
\mathcal{S}+\sum_{i=1}^n m_i\partial A_i'$ and $\gamma$ is less than
$2+2\log_2 K'$.  As $\gamma$ is isotopic to a component of $\partial
A_i'$ for some $i$, in the case that $m_1< K\sigma$, the distance
between $\partial \mathcal{S}+\sum_{i=1}^n m_i\partial A_i'$ and some
$\partial A_i'$ is bounded by a number that depends only on $K$,
$\partial\mathcal{S}$ and the $\partial A_i'$'s.

Therefore, combining the two cases above, \fullref{Lcurve} holds.
Moreover, it follows from the proof that the diameter of the set
$\{\partial\mathcal{S}+\sum_{i=1}^n m_i\partial A_i'\}$ can be found
algorithmically.

Now \fullref{T02} follows from \fullref{Lcurve} and the
discussions in \fullref{S0-eff}.

\begin{theorem1}
Let $M$ be a simple 3--manifold with connected boundary and a
0--efficient triangulation.  Let $S_k$ be the set of normal and almost
normal surfaces satisfying the following two conditions
\begin{enumerate} 
\item the boundary of each surface in $S_k$ consists of essential
curves in $\partial M$
\item the Euler characteristic of each surface in $S_k$ is at least
$-k$.
\end{enumerate} 
Let $C_k$ be the set of boundary curves of surfaces in $S_k$.  Then
$C_k$ has bounded diameter in the curve complex of $\partial M$.
Moreover, there is an algorithm to find the diameter.
\end{theorem1}
\begin{proof}
Let $S$ be a normal or an almost normal surface with $-\chi(S)\le k$.
So we have $S=\mathcal{S}+C+\sum m_iA_i$, where $C$ is a closed
surface and $A_i$ is a normal annulus in the fundamental solution.
Moreover, by \fullref{Pbounded}, there are only finitely many
possible surfaces for $\mathcal{S}$.

If we fix a $\mathcal{S}$, then \fullref{Lcurve} says that
$\{\partial S=\partial\mathcal{S}+\sum m_i\partial A_i\}$ has bounded
diameter.  Since there are only finitely many choices for
$\mathcal{S}$, $C_k$ has bounded diameter.  It follows from the proof
that there is an algorithm to find this diameter.
\end{proof}

\begin{proof}[Proof of \fullref{Tmain}]
\fullref{Tmain} follows immediately from \fullref{T02} and the
discussions in \fullref{Sstr} and \fullref{Sinter}.  By
\fullref{CStr}, there is surface $S_i$ ($i=1,2$) properly
embedded in $M_i$, such that $S_i$ is either essential or
$\partial$--strongly irreducible in $M_i$ and the distance
$d_{\mathcal{C}(F)}(\phi(\partial S_1), S_2)$ is at most $2g-2$, where
$g=g(M_1)+g(M_2)-g(F)$.  By \cite{FO} and a theorem in \cite{B} (see
the Appendix below for a workaround for \cite{B}), $S_i$ is isotopic
to a normal or an almost normal surface for any 0--efficient
triangulation of $M_i$, see \fullref{R1}.  Now we choose a
0--efficient triangulation for $M_i$ and \fullref{Tmain} follows
from \fullref{T02}.
\end{proof}

\section*{Appendix}

The purpose of this appendix is to address an issue in the proof of
\cite[Corollary 8.9]{B}.  While Bachman insists the proof is correct,
there is a concern on the thin-position argument for manifolds with
boundary in the proof of \cite[Lemma 8.5]{B}.  The following is a
workaround suggested by the referee.

Note that an essential surface is isotopic to a normal surface with
respect to any triangulation, so the issue here is on
$\partial$--strongly irreducible surfaces.  Suppose $S_1$ is a
strongly irreducible and $\partial$--strongly irreducible surface
properly embedded in $M_1$ as in \fullref{LBSS} and
\fullref{CStr}.  It follows from the sweepout argument in
\cite{BSS} that $S_1$ is compressible on both sides (see
\fullref{Dst}), and $\partial S_1$ consists of essential curves
in $\partial M_1$ (see the proof of \fullref{CBSS}).  Next we
show that $S_1$ does not admit nested $\partial$--compressions.

Suppose $S_1$ admits nested $\partial$--compressions, then we can find
a disk $D$ such that $\partial D=\alpha\cup\beta$,
$\alpha\subset\partial M_1$, $\beta\subset S_1$, and $\int(D)\cap
S_1\ne\emptyset$ consists of non-nested arcs.  Let
$\beta_1,\dots,\beta_k$ be the arcs of $\int(D)\cap S_1$ and $\delta_i$
($i=1,\dots,k$) the subdisk of $D$ bounded by $\beta_i$ and a subarc
of $\alpha$. By our assumption, $\delta_i\cap\delta_j=\emptyset$ if
$i\ne j$.  We may assume that each $\delta_i$ is a
$\partial$--compressing disk on the same side of $S_1$.  Moreover, we
may choose $D$ so that $k>0$ and $k$ is minimal among all such disks
$D$.  Let $Q=D-\bigcup_{i=1}^{k}\int(\delta_i)$.  Since $S_1$ is
compressible on both sides, there is a compressing disk $D'$ on the
opposite side of $\delta_i$ or equivalently on the same side as $Q$.
We may assume $D'\cap Q$ contains no closed curve. Since $S_1$ is
$\partial$--strongly irreducible, $D'\cap\beta_i\ne\emptyset$ for each
$i$.  Let $\gamma$ be an arc of $D'\cap Q$ that is outermost in $D'$
and $\Delta$ the subdisk of $D'$ cut off by $\gamma$ with $\Delta\cap
Q=\gamma$.  The arc $\gamma$ cuts $Q$ into two disks $Q_1$ and $Q_2$.
Thus either (1) $Q_i\cup\Delta$ ($i=1$ or $2$) is a compressing disk
disjoint from some $\delta_i$, a contradiction to the
$\partial$--strong irreducibility, or (2) the union of $Q_i\cup\Delta$
($i=1$ or $2$) and some $\delta_j$'s form a new disk similar to $D$,
which contradicts the assumption that $k$ is minimal.  Thus $S_1$ does
not admit nested $\partial$--compressions.

We call the two sides of $S_1$ plus and minus sides.  By the
definition of strongly irreducible surfaces (\fullref{Dst}),
$S_1$ is compressible on both sides.  If we perform a maximal
compression on the plus side of $S_1$ and discard the closed surface
components, then we get a surface $S_1^+$.  Since $S_1$ is
$\partial$--strongly irreducible, $S_1^+$ is incompressible and
$\partial$--incompressible on the minus side.  This basically follows
from \cite{CG}, see part (1) of \cite[Lemma 5.5]{S1} for a proof for
surfaces with boundary.  Note that the proof of part (1) of
\cite[Lemma 5.5]{S1} does not mention $\partial$--compressing disks
because the surface in \cite{S1} is strongly irreducible but may not
be $\partial$--strongly irreducible.  However, with the assumption of
$\partial$--strong irreducibility, the same proof of \cite{S1} shows
$S_1^+$ is also $\partial$--incompressible on the minus side.  Thus
either $S_1^+$ consists of $\partial$--parallel surfaces, or after
some $\partial$--compressions on the plus side, $S_1^+$ becomes an
essential surface $S_1'$ in $M_1$ with $d(\partial S_1^+,\partial
S_1')\le -\chi(S_1^+)$.  As the argument for essential surfaces does
not use Bachman's theorem \cite{B}, we may assume $S_1^+$ is
$\partial$--parallel.  Since $S_1$ has no nested
$\partial$--compressions, the $\partial$--parallel components of
$S_1^+$ are not nested.  We can also apply the same argument on the
minus side of $S_1$.  Therefore we may assume that $S_1$ is a
boundary-Heegaard surface as in \cite{B}.

Next we explain the controversial part of \cite{B} which is pointed
out by the referee.  The proof of the main theorem in \cite{B} is
basically a thin-position argument in which the 1--skeleton of the
triangulation is in thin position with respect a sweepout
$\{\mathcal{S}_t\}$ of a boundary-Heegaard surface.  A problem arises
when a thick level surface admits a high disk $D'$ and a low disk $D$
with $D\cap D'=\emptyset$.  If both $D$ and $D'$ lie in the interior
of $M_1$ then a simple isotopy as in \cite[Figure 4]{B} can reduce the
width. The controversial part in \cite{B} is the case that $D'$ lies
in the interior of $M_1$ and $D$ has a boundary arc in $\partial M_1$.
Note that one can assume that $D$ is a $\partial$--compressing disk,
since otherwise there is a low disk totally in $\partial M_1$ disjoint
from $D'$ and the usual isotopy can reduce the width.  However, if $D$
is a $\partial$--compressing disk, the usual width-reduction operation
as above would be a $\partial$--compression along $D$, which is not an
isotopy on the level surface any more.  Below is a workaround for this
situation.

We first glue a product $F\times I$ ($F=\partial M_1$) to $M_1$ and
obtain a manifold $M_1'$ ($M_1'\cong M_1$).  We can extend $S_1$ to a
surface $S_1'$ properly embedded in $M_1'$ by adding vertical annuli
in $F\times I$ along $\partial S_1$.  We fix a 0--efficient
triangulation of $M_1$ and suppose $F\times I$ is not triangulated.

We consider a special sweepout or foliation $\{\mathcal{S}_t\}$ as in
\cite{B} with the restriction that for each regular leaf
$\mathcal{S}_t$, $\mathcal{S}_t\cap(F\times I)$ is obtained by pushing
pairwise disjoint $\partial$--compressing disks of $S_1$ (on the same
side) into $F\times I$, and $\mathcal{S}_t\cap M_1$ is obtained by
$\partial$--compressions on one side.  Note that if a
$\partial$--compression on $S_1$ yields a $\partial$--parallel disk
component, we also push the disk component into $F\times I$.

We now apply the thin-position argument on $\{\mathcal{S}_t\cap M_1\}$
and assume the 1--skeleton is in thin position.  Suppose a thick level
$\mathcal{S}_t$ admits a pair of disjoint high and low disks in the
2--skeleton.  Let $D$ be the low disk as explained above and suppose
$\partial D=\alpha\cup\beta$ with $\alpha\subset\mathcal{S}_t$ and
$\beta\subset\partial M_1$.  We may assume $D$ is a
$\partial$--compressing disk for $\mathcal{S}_t\cap M_1$.  Since $S_1$
is $\partial$--strongly irreducible, the high disk lies in the
interior of $M_1$.  Thus we can perform an isotopy on the
triangulation as in \cite[Figure 4]{B} by pushing the high disk down
and the low disk up, which leads to a contradiction to the
thin-position assumption.  Note that the isotopy of pushing the low
disk $D$ into $F\times I$ can be viewed as a $\partial$--compression
on $\mathcal{S}_t$ in $M_1$.  Moreover, by the assumptions on $S_1$,
after a $\partial$--compression on one side, there is no
$\partial$--compressing disk in $M_1$ on the other side, hence all the
$\partial$--compressions are on the same side.

Therefore the arguments in \cite{B,Sto} imply that one can isotope
$S_1'$ into a surface $\Sigma'$ so that (1) $\Sigma'\cap (F\times I)$
is obtained by pushing pairwise disjoint $\partial$--compressing disks
of $S_1$ (on the same side) into $F\times I$, and $\Sigma=\Sigma'\cap
M_1$ is obtained by $\partial$--compressions on one side, note that we
also push the possible trivial disk components into $F\times I$, and
(2) $\Sigma$ is normal or almost normal with respect to the
triangulation of $M_1$.  Furthermore, after one more
$\partial$--compression in a tetrahedron, we may assume the special
type of almost normal pieces in Figure 9 of \cite{B} does not appear.
This implies that $\partial\Sigma$ is normal in $\partial M_1$.

Now we study the property of $\Sigma$ and $\Sigma'\cap(F\times I)$.
First, since all the $\partial$--compressions occur on the same side
of $S_1$, the trivial-circle components of $\partial\Sigma$ (if any)
are not nested.  Since $\partial\Sigma$ is normal and by part (a) of
\fullref{P1}, this implies that $\partial\Sigma$ contains at
most one trivial-curve component.  Note that
$\partial\Sigma\ne\emptyset$ since $F$ is incompressible and $S_1$ is
a subsurface of a Heegaard surface.

If $\partial\Sigma$ contains at least one essential curve, then as
before, the distance $d(\partial\Sigma,\partial S_1)<-\chi(S_1)$,
viewed in the curve complex $\mathcal{C}(F)$.  Now we can prove the
main theorem by applying the arguments in Sections \ref{S0-eff},
\ref{Sannuli} and \ref{Scurve} on $\Sigma$, see \fullref{R1} and
the remark before \fullref{P1}.

Therefore we may suppose $\partial\Sigma$ is a single trivial
vertex-linking curve in $\partial M_1$. Let $\delta$ be the disk
bounded by $\partial\Sigma$ in $\partial M_1$ and
$P=\Sigma'\cap(F\times I)$.  So $P\cup\Sigma=\Sigma'$ and by the
construction of $\Sigma$, $\delta\cup P$ is $\partial$--parallel in
$F\times I$, in other words, $P=\Sigma'-\int(\Sigma)$ can be
constructed by adding a vertical tube to a (once-punctured)
$\partial$--parallel surface in $F\times I$.

Let $F_-=\partial M_1-\int(\delta)$.  There is a natural projection
from the arc-and-curve complex $\mathcal{AC}(F_-)$ to the curve
complex $\mathcal{C}(\partial M_1)=\mathcal{C}(F)$ denoted by $\pi
\co \mathcal{AC}(F_-)\to\mathcal{C}(F)$ as follows.  We view
$F_-=F-\int(\delta)$.  For any closed essential curve $\gamma$ in
$F_-$, $\gamma$ is also an essential curve in $F$, we set $\pi
([\gamma])=[\gamma]$.  For any essential arc $\alpha$ in $F_-$, let
$\hat{\alpha}$ be the closed curve obtained by connecting
$\partial\alpha$ by an arc properly embedded in the disk $\delta$.  We
define $\pi ([\alpha])=[\hat{\alpha}]$.  Note that if two arcs
$\alpha\cap\beta=\emptyset$ in $F_-$, then
$\hat{\alpha}\cap\hat{\beta}$ is either empty or a single point.  This
means that if $d(\alpha,\beta)=1$ in $\mathcal{AC}(F_-)$ then
$d(\pi(\alpha),\pi(\beta))=d(\hat{\alpha},\hat{\beta})\le 2$ in
$\mathcal{C}(F)$.

Note that the disk $\delta$ is a compressing disk for $S_1'=\Sigma'$.
We denote the two sides of $S_1'$ using plus and minus and suppose
$\delta$ is on the plus side.  Since $S_1'$ is compressible on both
sides, there is another compressing disk $D$ on the minus side and
$\partial\delta\cap\partial D\ne\emptyset$ in $S_1'$.  Since $\Sigma$
is obtained by $\partial$--compressions on the plus side, $\Sigma$ is
$\partial$--incompressible on the minus side and $D\cap\partial
M_1\ne\emptyset$.  Moreover, $\partial M_1$ cuts $D$ into a collection
of subdisks and all the bigon disks lie in $F\times I$ (since $\Sigma$
is $\partial$--incompressible on the minus side in $M_1$).  Let $D_1$
be such a bigon subdisk of $D$ and suppose $\partial
D_1=\alpha\cup\beta$, where $\alpha\subset\partial M_1$ and
$\beta\subset P$.  Let $D_0$ be the subdisk of $D$ adjacent to $D_1$
with $D_0\cap D_1=\alpha$ and $D_0\subset M_1$.  Since $P$ can be
obtained by adding a vertical tube to a punctured $\partial$--parallel
surface in $F\times I$, $\hat{\alpha}$ and $\partial S_1'$ project to
disjoint curves in $F$, ie, $d(\pi (\alpha),\partial S_1')\le 1$.

Note that $\Sigma$ cuts $M_1$ into two submanifolds and we denote the
one on the minus side by $N$.  Clearly $D_0$ is a compressing disk for
$N$.  For any compressing disk $\Delta$ of $N$, since $\Sigma$ is
incompressible on the minus side,
$\partial\Delta\cap\partial\delta\ne\emptyset$, ie, $\partial\delta$
is disk-busting in $N$.  For any compressing disk $\Delta$ of $N$, we
suppose $|\partial\Delta\cap\partial\delta|$ is minimal among all
disks in the isotopy class of $\Delta$.  We fix an arc component
$\gamma_\Delta$ of $\partial\Delta\cap\partial M_1$ for each $\Delta$.
Let $\mathcal{D}$ be the disk complex of $\partial N$ (ie, curves of
$\partial N$ bounding compressing disks in $N$). Define a projection
$\pi_A\co \mathcal{D}\to\mathcal{AC}(F_-)$ as
$\pi_A([\partial\Delta])=[\gamma_\Delta]$.  The following theorem in
\cite{L3} was also independently proved by Masur and Schleimer.

\begin{theorem*}[\cite{L3}]
 Let $N$ be as above, $\mathcal{D}$ the disk complex, and $F_-$ a
 compact essential subsurface of $\partial N$. Suppose $\partial F_-$
 is disk-busting in $\partial N$.  Then either
\begin{enumerate}
\item $N$ is an $I$--bundle of which $F_-$ is a horizontal boundary
component, or
\item the image $\pi_A(\mathcal{D})$ of the disk complex has diameter
at most 10 in $\mathcal{AC}(F_-)$ and $\pi\circ\pi_A(\mathcal{D})$ has
diameter at most 20 in $\mathcal{C}(F)$.
\end{enumerate}
\end{theorem*}

Note that part (a) of the theorem cannot happen in our case because
otherwise one could isotope $F$ to be disjoint from the Heegaard
surface.  Thus for any compressing disk $\Delta$ of $N$,
$d(\hat{\gamma}_\Delta,\hat{\alpha})\le 20$ in $\mathcal{C}(F)$, where
$\alpha$ is the arc $D_1\cap D_0$ above.  Moreover, since
$d(\hat{\alpha},\partial S_1)\le 1$, we have
$d(\hat{\gamma}_\Delta,\partial S_1)\le 21$ for any compressing disk
$\Delta$ of $N$.

Let $\Gamma$ be the set of almost normal surfaces in $M_1$ such that
for each surface $X$ in $\Gamma$, $\partial X$ is a vertex linking
circle in $\partial M_1$ and $\chi(X)\ge\chi(S_1)$.  As in
\fullref{S0-eff} and \fullref{Sannuli}, there is a finite
collection of branched surfaces such that each surface in $\Gamma$ is
fully carried by a branched surface in the collection, and for each
branched surface $B$, $\partial B$ a single trivial circle in
$\partial M_1$.  For any surface $X$ in $\Gamma$, $\partial X$ bounds
a disk $\delta$ in $\partial M_1$.  Let $N$ be the closure of the
component of $M_1-X$ that contains $F_-=\partial M_1-\delta$.  If $X$
is fully carried by $B$, then $N$ can be constructed by connecting
some components of $M_1-\int(N(B))$ using $I$--bundles.  Although there
may be infinitely many surfaces in $\Gamma$, since there are only
finitely many branched surfaces and $\chi(X)$ is bounded, there are
only finitely many possible topological types for $N$ and we can list
them all.  For each possible $N$, we randomly find a compressing disk
$\Delta$ for $N$ and fix an arc $\gamma_\Delta$ of $\partial\Delta\cap
F_-$.  So we can construct finitely many closed curves
$\hat{\gamma}_\Delta$.  By the discussion above, if the gluing map
$\phi\co \partial M_1\to\partial M_2$ is so complex that
$d(\hat{\gamma}_\Delta,\partial S_1)> 21$ for each possible $N$, then
it is impossible to have a surface $S_1$ with all the
requirements. This implies that the original Heegaard surface cannot
be strongly irreducible and the theorem follows.

\bibliographystyle{gtart} \bibliography{link}

\end{document}